\documentclass{amsart}

\usepackage{amsfonts,amsmath}
\usepackage{amssymb,latexsym}

\newtheorem{theorem}{Theorem}%[section]
\newtheorem{lemma}[theorem]{Lemma}%[section]
\newtheorem{definition}[theorem]{Definition}%[section]
\newtheorem{proposition}[theorem]{Proposition}%[section]
\newtheorem{corollary}[theorem]{Corollary}%[section]
\newtheorem{remark}[theorem]{Remark}%[section]
\numberwithin{equation}{section} 
\newcommand{\real}{\mathbb{R}}
\newcommand{\loc}{\scriptstyle{loc}} 
\newcommand{\vare}{\varepsilon}

\begin{document} 
\title[2D incompressible viscous flow around a small obstacle]{Two-dimensional incompressible viscous flow around a small obstacle} 
\author[Iftimie, Lopes Filho and Nussenzveig Lopes]{{D. Iftimie,}
{M. C. Lopes Filho and}  
{H. J. Nussenzveig Lopes}}

\begin{abstract}
In this work we study the asymptotic behavior of viscous incompressible 2D flow in the
exterior of a small material obstacle. We fix the initial vorticity $\omega_0$ and the circulation
$\gamma$ of the initial flow around the obstacle. We prove that, if $\gamma$ is sufficiently small, the 
limit flow satisfies the full-plane Navier-Stokes system, with initial vorticity $\omega_0 + 
\gamma \delta$, where $\delta$ is the standard Dirac measure. The result should be contrasted
with the corresponding inviscid result obtained by the authors in \cite{iln-shrink}, where the
effect of the small obstacle appears in the coefficients of the PDE
and not only on the initial data.
The main ingredients of the proof are $L^p-L^q$ estimates for the Stokes operator in an exterior
domain, a priori estimates inspired on Kato's fixed point method, energy estimates, renormalization 
and interpolation.
\end{abstract}

\maketitle

\tableofcontents

\section{ Introduction}

The purpose of this work is to study the influence of a material obstacle on
the behavior of two-dimensional incompressible viscous flows when the size of the 
obstacle is small compared to that of a reference spatial scale. More precisely, we fix
both an initial vorticity $\omega_0$, smooth and compactly supported, and the circulation
$\gamma$ of the initial velocity around the boundary of the obstacle, while 
homothetically contracting the obstacle to a point $P$ outside the support of $\omega_0$. 
The initial vorticity $\omega_0$ and the circulation $\gamma$ uniquely determine a 
family of divergence-free initial velocities $u_0^{\vare}$ with 
$\mbox{curl }u^{\vare}_0 = \omega_0$ and $u^{\vare}_0(x) \to 0$ at infinity. The size 
of the support of the initial vorticity $\omega_0$ can be used as reference spatial scale. Let 
$u^{\vare} = u^{\vare}(x,t)$ be a solution of the Navier-Stokes equations with
initial data $u_0^{\vare}$ and no-slip data at the boundary of the small obstacle. 
Our problem is to determine the asymptotic behavior of $u^{\vare}$ as $\vare \to 0$.
We will show that $u^{\vare}$ converges to a solution of the Navier-Stokes equations in the 
full plane with initial vorticity $\omega_0 + \gamma \delta(x-P)$, as long as 
$\gamma$ is sufficiently small.

There is a sharp contrast between the behavior of ideal and viscous flows around a 
small obstacle. In \cite{iln-shrink}, the authors studied the vanishing obstacle problem for 
incompressible, ideal, two-dimensional flow. The ideal flow assumption is physically incorrect 
in the presence of material boundaries, and part of the motivation of the present work,
together with \cite{iln-shrink}, is to explore more precisely the extent to which the ideal
flow assumption misrepresents the physical flow in the presence of material boundaries.  
The main result in \cite{iln-shrink} is that the  limit vorticity in the ideal case satisfies a modified 
vorticity equation of the form $\omega_t + u\cdot\nabla\omega=0$, with $\mbox{div }u = 0$ 
and $\mbox{curl }u = \omega + \gamma \delta(x-P)$. In the ideal case, the correction due to the obstacle
appears as time-independent additional convection centered at $P$, whereas in the viscous case, the
correction appears only on the initial data and gets convected and diffused.   

The small obstacle limit is an instance of the general problem of PDE on singularly 
perturbed domains. There is a large literature on such problems, specially in
the elliptic case, see \cite{MNP00} for a broad overview. Asymptotic behavior
of fluid flow on singularly perturbed domains is a natural subject for analytical
investigation which is virtually unexplored. The present work, together with 
\cite{iln-shrink}, may be regarded as a first attempt to address this class of problems.

There is a natural connection between the approximation problem as we have formulated it 
and the issue of uniqueness for the limit problem. In fact, from a technical point of view,
our work is closely related to the classical uniqueness result due to Y. Giga, T. Miyakawa and 
H. Osada, on solutions of the incompressible 2D Navier-Stokes equations with measures as 
initial data, see \cite{GMO}. Some of the more striking similarities are: the difficulties 
with locally infinite kinetic energy, the use of $L^p$ estimates  for the linearized problem and
the use of Kato-type norms to estimate the nonlinearity. The smallness condition on the mass 
of the point vortices in the initial data, required in the uniqueness result, is closely related
to our smallness condition on the circulation. 

The remainder of this work is organized in ten sections. In Section 2 we summarize $L^p$ estimates for the
time-dependent Stokes problem on exterior domains. In Section 3 we formulate precisely the problem we
wish to discuss and write uniform estimates for the initial data. In Section 4 we study the asymptotic
behavior of the initial data. In Section 5 we discuss physical motivation for our problem and we establish
the small obstacle asymtotics for circularly symmetric flows, a linear version of our problem. In Section 
6 we derive a priori estimates in the initial layer for the nonlinear correction term. In Section 7 we 
deduce global-in-time energy estimates for the nonlinear correction term. In Section 8 we put together 
the estimates for the linear part with the estimates for the nonlinear correction, obtaining a complete
set of a priori estimates for velocity. In Section 9 we prove compactness space-time, in Section 10
we perform the passage to the limit and in Section 11 we add comments and concluding remarks.

We conclude this introduction with a few remarks regarding notation. Given a vector 
$z=(z_1,z_2) \in \real^2$ we denote its orthogonal vector by $z^{\perp} = (-z_2,z_1)$. We use the subscript 
$c$ in function spaces to denote compact support, as in $C^{\infty}_c$, and we use standard notation for 
Sobolev spaces, $W^{k,p}$, where $1 \leq p \leq \infty$ and $k \in \mathbb{Z}$, with $H^{k}$ standing for 
the case $p=2$. We use the subscript $loc$ in function spaces $X$ to denote functions which are locally 
in $X$. In particular, $L^p_{\loc}([0,\infty);W^{k,q})$ denotes functions $f = f(t,x) \in L^p([0,M];W^{k,q})$ 
for any $M>0$, whereas $L^p_{\loc}((0,\infty);W^{k,q})$ denotes functions $f = f(t,x) \in L^p([\delta,M];W^{k,q})$ 
for any $\delta>0$ and any $M>0$, but not necessarily for $\delta = 0$. Finally, $L^{2,\infty}$ denotes the 
Lorentz space of functions $f$ whose distribution function satisfies 
$\lambda_f=\lambda_f(s)=|\{|f|>s\}| = \mathcal{O}(s^{-2})$.
 
\section{ Estimates for the Stokes semigroup}

    In this section we will put together several results on estimates for the
Stokes semigroup on exterior domains. Let us begin by introducing some basic notation.

Let $\Omega$ be a bounded, open, simply connected subset of $\mathbb{R}^2$ with boundary $\Gamma$,
a smooth Jordan curve. We denote by $\Pi$ the unbounded connected component of $\mathbb{R}^2\setminus\Gamma$.
Fix $\nu >0$ and let $\mathbb{P}$ denote the Leray projector onto divergence-free vector fields on $\Pi$. 
Let $\mathbb{A} \equiv -\mathbb{P} \Delta$ be the Stokes operator on $\Pi$ and denote the Stokes semigroup by
$S_{\nu}(t) = e^{-\nu t\mathbb{A}}$. Given $v_0 \in C^{\infty}_c(\Pi)$, let $v(t,x) = S_{\nu}(t) v_0$ be the
unique solution of the system
\begin{equation} \label{stokes}
\left\{
\begin{array}{ll}
\partial_t v - \nu \Delta v = - \nabla p, & \mbox{ in } (0,\infty) \times \Pi  \\
\mbox{div }v = 0,& \mbox{ in } [0,\infty) \times \Pi   \\
v = 0,& \mbox{ on } (0,\infty) \times \Gamma   \\
\lim_{|x| \to \infty} \, v(t,x) = 0, & \mbox{ for all } t \geq 0\\
v(0,x) = v_0(x),& \mbox{ on } \{t=0\} \times \Pi. 
\end{array}   
\right.
\end{equation}  

We denote by $X^p(\Pi)$ the closure of the space of divergence-free, $C^{\infty}_c(\Pi)$ 
vector fields with respect to the $L^p$-norm. The Stokes operator in $X^p$ 
generates an analytic semigroup of class $C^0$ on $X^p(\Pi)$, for any $1<p<\infty$, see \cite{giga81},
so that, in particular, problem \eqref{stokes} is well-posed in $X^p(\Pi)$. 

We will require two kinds of estimates on the 
Stokes semigroup, $L^p$ estimates and renormalized energy estimates. We first state the $L^p$ estimates. 

\begin{theorem} \label{ellpest}
Let $1 < q <\infty$. Consider $v_0 \in X^q(\Pi)$ and $F \in L^q(\Pi;M_{2\times2}(\real))$. 
Then we have the following estimates.

\begin{enumerate}
\item[(S1)] Let $q \leq p<\infty$. There exists $K_1 = K_1(\Pi,p,q)>0$ such that
\[ \|S_{\nu}(t)v_0\|_{L^p} \leq K_1 (\nu t)^{\frac{1}{p} - \frac{1}{q}}\|v_0\|_{L^q},\] 
for all $t>0$.
\item[(S2)] Let $q \leq p \leq 2$. There exists $K_2 = K_2(\Pi,p,q)>0$ such that
\[ \|\nabla S_{\nu}(t)v_0\|_{L^p} \leq K_2 (\nu t)^{-\frac{1}{2} + \frac{1}{p} - \frac{1}{q}}\|v_0\|_{L^q},\]
for all $t>0$.
\item[(S3)] Assume $q \geq 2$ and let $q \leq p < \infty$. Then there exists $K_3 = K_3(\Pi,p,q)>0$
such that
\[ \| S_{\nu}(t)\,\mathbb{P} \mbox{ div }F\|_{L^p} \leq K_3 (\nu t)^{-\frac{1}{2} + \frac{1}{p} - \frac{1}{q}}
\|F\|_{L^q},\]
for all $t>0$, with the divergence taken along rows of the matrix $F$.
\end{enumerate} 
\end{theorem} 

This theorem summarizes several results already contained in the literature, which we
have collated above for convenience. 

\begin{proof}
Estimates (S1) and (S2) were proved
in \cite{DS1,MS} (see also \cite{DS2} for the case $p = \infty$). Estimate (S3) follows from (S2) by duality. Indeed, the adjoint of 
$S_{\nu}(t)$ on $X^p$ is again $S_{\nu}(t)$, defined on $X^{p'}$, with $1/p + 1/p' = 1$ and 
therefore the adjoint of $\nabla S_{\nu}(t)$ is $S_{\nu}(t)\, \mathbb{P} \mbox{ div }$.
The dependence on the viscosity follows directly by rescaling time, 
since $S_{\nu}(t)= S_1(\nu t)$ .
\end{proof}

Next we address a renormalized energy estimate for the Stokes semigroup. Our concerns include 
infinite energy solutions to the Navier-Stokes equations whose 
behavior at infinity is $\mathcal{O}(1/|x|)$. In the following result we will prove that 
solutions to the Stokes system retain the behavior at infinity of their initial data.  

\begin{proposition} \label{renen}
Let $v_0$ be a smooth divergence-free vector field on $\Pi$ vanishing at the boundary $\Gamma$. 
We assume also that $v_0 \in X^p(\Pi)$  for some $p>2$ and that $\nabla v_0 \in L^2(\Pi)$. Then
\[S_{\nu}(t)v_0 - v_0 \in C^0([0,\infty);L^2(\Pi)) \cap L^2_{\loc}([0,\infty);H^1(\Pi)).\]
Moreover the following inequality holds
\begin{equation} \label{renendesig}
\|S_{\nu}(t)v_0 - v_0\|_{L^2}^2 + \nu \int_0^{t} \|\nabla [S_{\nu}(\tau)v_0 - v_0]\|_{L^2}^2 d\tau \leq \nu t 
\|\nabla v_0\|_{L^2}^2. 
\end{equation}
\end{proposition}  

\begin{proof}
Let $W = S_{\nu}(t)v_0 - v_0$. Then $W$ satisfies:
\begin{equation} \label{grac} \left\{
\begin{array}{ll}
\partial_t W - \nu \Delta W = - \nabla p + \nu \Delta v_0, & \mbox{ in } (0,\infty) \times \Pi  \\
\mbox{div }W = 0,& \mbox{ in } [0,\infty) \times \Pi   \\
W = 0,& \mbox{ on } (0,\infty) \times \Gamma   \\
\lim_{|x| \to \infty} \, W(t,x) = 0, & \mbox{ for all } t \geq 0\\
W(0,x) = 0, & \mbox{ on } \{t=0\} \times \Pi. 
\end{array}   \right.\end{equation}

It is well-known that \eqref{grac} admits a unique solution 
$\widetilde{W}$ in $C^0([0,\infty);L^2(\Pi)) \cap L^2_{\loc}([0,\infty);H^1(\Pi))$, 
see, for instance, Theorem III.1.1 in \cite{temam84}. The fact that $W-\widetilde{W}=0$ 
follows from the well-posedness of \eqref{stokes} in $X^p$.  The standard energy 
estimate gives \eqref{renendesig}.
\end{proof}

One consequence of the nontrivial topology of $\Pi$ is the existence of harmonic
vector fields, i.e. divergence-free and curl-free vector fields which are tangent to
$\Gamma$ and vanish at infinity. We denote by $H_{\Pi}$ the unique harmonic vector field 
on the exterior domain $\Pi$ which satisfies the condition
\[ \oint_{\Gamma} H_{\Pi} \cdot ds = 1,\]
where the contour integral is taken in the counterclockwise sense.
It is an elementary application of Hodge theory that the  vector space
of these
harmonic 
vector fields on $\Pi$ is one dimensional, and we can use $H_{\Pi}$ as a basis.      
In the case where $\Pi$ is the exterior of the unit disk centered at the origin,
we will denote $H_{\Pi}$ simply by $H$, and we have: 
\begin{equation} \label{helena}
H = \frac{x^{\perp}}{2\pi|x|^2}.
\end{equation}
 
We will require detailed information on the behavior of $H_{\Pi}$ both at infinity and
near the boundary $\Gamma$, which we obtain by means of a conformal mapping.  We denote 
\[\mathcal{U}  \equiv \{|x|>1\}\] 
and switch to complex variables notation in the result below.  

\begin{lemma} \label{confmap} There exists a smooth biholomorphism $T:\Pi \rightarrow \mathcal{U}$, 
extending smoothly up to the boundary, mapping $\Gamma$ to $\{|z|=1\}$. Furthermore, there 
exists a nonzero real number $\beta$ and a bounded holomorphic function $h: \Pi \to \mathbb{C}$ 
such that:
\begin{equation} \label{holomform}
T(z) = \beta z + h(z).
\end{equation}
Additionally, 
\begin{equation} \label{morestuff}
h^{\prime}(z) = {\mathcal O}\left(\frac{1}{|z|^2}\right), \mbox{ as } |z| \to \infty.
\end{equation}
\end{lemma}

This Lemma is an excerpt from \cite{iln-shrink}. Its proof is an exercise in complex 
analysis.

It was observed in \cite{iln-shrink} (see identity (2.10) in \cite{iln-shrink}) that
\begin{equation} \label{harm}
H_{\Pi} = H_{\Pi}(x) = \frac{1}{2\pi} \frac{DT^t(x) (T(x))^{\perp}}{|T(x)|^2}.
\end{equation}
From Lemma \ref{confmap}, we see that $|H_{\Pi}|$ is $\mathcal{O}(1/|x|)$ for large $|x|$.
This implies that $H_{\Pi}$ belongs to the Lorentz space $L^{2,\infty}(\Pi)$.

We close this section with an estimate for the Stokes semigroup acting on infinite energy 
initial data. 

\begin{proposition} \label{stokesagapi}
Let $2 < p < \infty$ and let $v_0 \in L^{2,\infty}(\Pi) \cap X^p(\Pi)$. 
There exists a constant $K_5>0$ such that 
\[\|S_{\nu}(t)v_0\|_{L^p} \leq K_5 (\nu t)^{\frac{1}{p}-\frac{1}{2}} \|v_0\|_{L^{2,\infty}}.\]
In particular, this estimate holds true for $v_0 = H_{\Pi}(x)$.
\end{proposition}

\begin{proof}
This estimate is contained in Proposition 2.2, item (4), of \cite{kozono-yamazaki}. 
To see that it holds for $H_{\Pi}$, we first show that $H_{\Pi} \in X^p(\Pi)$ for any $p>2$.
This is easy to prove in the case $\Pi = \mathcal{U}$ because, for any 
function $\varphi \in C^{\infty}_c((0,\infty))$, $\varphi(|x|) H(x)$ is smooth,
compactly supported and divergence-free, and, by taking $\varphi^{\vare}$ 
a sequence of cutoffs for the interval $(1+\vare,1/\vare)$, it
is easy to see that $\varphi^{\vare}(|x|) H \to H$ in $L^p$, for $p>2$.
For general $\Pi$, we use the conformal mapping $T$, approximating 
$H_{\Pi}$ by $\varphi^{\vare}(|T(x)|)H_{\Pi}(x)$, where $\varphi^{\vare}$
is the same family of cutoffs used in the case of the exterior of the disk.
This strategy works because $\varphi^{\vare}(|T(x)|)H_{\Pi}(x)$ is also divergence-free.
\end{proof}

\section{ The evanescent obstacle}

The purpose of this section is to set down a precise statement of the small
obstacle problem. Many of the key issues 
regarding the small obstacle limit and incompressible flow have been 
discussed in detail in \cite{iln-shrink}, so that we will focus on
issues specifically related with viscous flow and briefly outline the rest. 

As in \cite{iln-shrink}, fix $\omega_0 \in C^{\infty}_c(\real^2)$
and assume that the origin does not belong to the support of $\omega_0$. Let $\Omega$ be
a bounded, open, connected and simply-connected subset of the plane whose boundary $\Gamma$
is a $C^{\infty}$ Jordan curve. The evanescent obstacle is the family of domains 
$\vare \Omega$, with $0 < \vare < \vare_0$. The parameter $\vare_0$ 
is chosen small enough so that the support of $\omega_0$ does not intercept 
$\vare\Omega$ for any $0 < \vare < \vare_0$.

Fix $0 < \vare < \vare_0$.
Let $\Pi_{\vare} \equiv \real^2 \setminus \overline{\vare\Omega}$ and 
$\Gamma_{\vare} = \partial \Pi_{\vare}$. We use the conformal mapping
$T : \Pi_1 \to \mathcal{U}$, given in Lemma \ref{confmap}, to define a family of smooth
biholomorphisms 
\begin{equation} \label{tvare}
T^{\vare} = T^{\vare}(x) \equiv T\left(\frac{x}{\vare}\right). 
\end{equation}
Throughout we write  $H^{\vare}$ for $H_{\Pi_{\vare}}$ and 
$G^{\vare} = G^{\vare}(x,y)$ will be the Green's function of the Laplacian in $\Pi_{\vare}$.
Let $K^{\vare}(x,y) = \nabla^{\perp}_x G^{\vare}(x,y)$ be the kernel 
of the Biot-Savart law on $\Pi_{\vare}$ and denote the associated integral operator 
by $f \mapsto K^{\vare}[f] = \int_{\Pi_{\vare}} K^{\vare}(x,y)f(y)\,dy$. 
Both $K^{\vare}$ and $H^{\vare}$ are related to $K_{\mathcal{U}}$ and $H_{\mathcal{U}}$ respectively, 
through the conformal mapping $T^{\vare}$, in a way which was made explicit in 
\cite{iln-shrink}. The relevant fact is the way that both the Biot-Savart kernel
and the basic harmonic vector field scale with $\vare$,  see identities (3.5) and (3.6) 
in \cite{iln-shrink}.

Fix $\alpha \in \real$ and let
\begin{equation} \label{u0}
u_0^{\vare} \equiv K^{\vare}[\omega_0] + \alpha H^{\vare}. 
\end{equation}

We consider the problem
\begin{equation} \label{NSevan}
\left\{
\begin{array}{ll}
\partial_t u^{\vare} + u^{\vare} \cdot \nabla u^{\vare} - \nu \Delta u^{\vare}= 
- \nabla p^{\vare}    , & \mbox{ in } (0,\infty) \times \Pi_{\vare}  \\
\mbox{div } u^{\vare} = 0,& \mbox{ in } [0,\infty) \times \Pi_{\vare}   \\
u^{\vare} = 0,& \mbox{ on } (0,\infty) \times \Gamma_{\vare}   \\
\lim_{|x| \to \infty} \, u^{\vare} (t,x) = 0, & \mbox{ for all } t \geq 0\\
u^{\vare}(0,x) = u^{\vare}_0(x),& \mbox{ on } \{t=0\} \times \Pi_{\vare}. 
\end{array}   
\right.
\end{equation}

We begin by observing that $u_0^{\vare} \in L^{2,\infty}(\Pi_{\vare}) \cap 
L^p(\Pi_{\vare})$ for any $2<p\leq \infty$. Indeed, $u_0^{\vare}$ is smooth, 
and therefore locally bounded, so that we only require knowledge on the behavior of 
$u_0^{\vare}$ at infinity. By Lemma \ref{confmap} and identity \eqref{harm} 
$|H^{\vare}|$ has $\mathcal{O}(1/|x|)$ behavior as $|x| \to \infty$, 
and therefore it belongs to $L^{2,\infty}(\Pi_{\vare}) \cap L^p(\Pi_{\vare})$ for 
any $2< p \leq \infty$. In fact, the $L^{2,\infty}$ bound on $|H^{\vare}|$ is independent
of $\vare$ as can be readily seen by rescaling to a fixed domain and
using that $H_{\Pi}$ belongs to $L^{2,\infty}$. 
In \cite{iln-shrink} it was shown that $|K^{\vare}[\omega_0]|$ 
has behavior $\mathcal{O}(1/|x|^2)$ at infinity (see estimate (2.8) in \cite{iln-shrink}) 
and therefore it belongs to $L^p(\Pi_{\vare})$, for any $p\geq 2$, 
and, in particular, to $L^{2,\infty}(\Pi_{\vare})$.

Global-in-time well-posedness for problem \eqref{NSevan} was established by Kozono 
and Yamazaki in \cite{kozono-yamazaki}. The existence part of Kozono and Yamazaki's
result requires that the initial velocity satisfy a smallness condition of the form
\[ \limsup_{R \to \infty} R \;|\{ x \in \Pi_{\vare} \;|\; |u_0^{\vare}(x)| > R\}|^{1/2} \ll 1.\]
Since $u_0^{\vare}$ is bounded, the limsup above is always zero, for any $\vare>0$. 
Uniqueness holds for divergence-free initial data in $L^{2,\infty} + X^p$ without any 
additional conditions.   

The evanescent obstacle problem consists of understanding the asymptotic 
behavior of Kozono and Yamazaki's solution $u^{\vare}(x,t)$ for small $\vare$. 
More precisely we will show that, under appropriate assumptions, $u^{\vare}$ has
a limit, and we will identify an equation satisfied by this limit.    

Fix $\varphi: \real \to [0,1]$ a smooth, monotone function such that $\varphi(s) \equiv 0$ if 
$s \leq 2$ and $\varphi(s) \equiv 1$ if $s \geq 3$. For each $\vare > 0$ and $\lambda >0$
we introduce the adapted cut-off functions:
\begin{equation} \label{cutoff}
\varphi^{\vare,\lambda}(x) \equiv \varphi\left(\frac{\vare}{\lambda} |T^{\vare}(x)|\right),
\end{equation}    
Note that the cutoff function $\varphi^{\vare,\lambda}$ vanishes in a ball of radius 
$\mathcal{O}(\lambda)$ and it is identically equal to $1$ outside a larger ball of radius 
$\mathcal{O}(\lambda)$, for large $\lambda$. Furthermore, the radii of the annulus 
where $\varphi^{\vare,\lambda}$ is not constant can be made independent of $\vare$. This follows
easily from the fact that $T$ is asymptotically affine at infinity, see \eqref{holomform}.

We will now introduce a pair of parameters that are useful to describe the asymptotic behavior
of $u^{\vare}_0$ when $\vare \to 0$. Consider 
\begin{equation} \label{mandgamma}
m \equiv \int_{\real^2} \omega_0 \, dx \hspace{1cm} \mbox{ and } \hspace{1cm} 
\gamma \equiv \oint_{\Gamma} u_0^{\vare} \cdot ds. 
\end{equation}
By Stokes' Theorem we have that $\gamma = \alpha - m$, and therefore, the 
circulation $\gamma$ does not the depend on $\vare$, see the proof of Lemma 3.1 in 
\cite{iln-shrink}.

For each $\lambda>0$, we introduce a convenient decomposition of the initial velocity as
\[u^{\vare}_0 = \mathbf{b}^{\vare}_0 + \mathbf{i}^{\vare}_0 + \mathbf{o}^{\vare}_0,\] 
with
\[\mathbf{b}_0^{\vare} \equiv K^{\vare}[\omega_0] + m (1 - \varphi^{\vare,\lambda}) H^{\vare}, \] 
\[ \mathbf{i}_0^{\vare} \equiv \gamma(1 - \varphi^{\vare,\lambda}) H^{\vare}, \] 
and 
\[ \mathbf{o}_0^{\vare} \equiv \alpha \varphi^{\vare,\lambda} H^{\vare}. \]
 
We need  to understand the behavior of each of the components of this decomposition, in the limit 
$\vare \to 0$. This is the content of our next result. The proof uses a large part of the work done 
in \cite{iln-shrink}.

\begin{lemma} \label{idp}
There exists $\lambda_0 > 0$, independent of $\vare$, for which $\|\mathbf{b}_0^{\vare}\|_{L^{p_1}}$,
$\|\mathbf{i}_0^{\vare}\|_{L^{p_2}}$ and $\|\mathbf{o}_0^{\vare}\|_{L^{p_3}}$ are uniformly
bounded in $\vare$, for any $1 < p_1 \leq \infty$, $1 \leq p_2 <2$ and 
$2 < p_3 \leq \infty$. The vector fields $\mathbf{b}_0^{\vare}$, $\mathbf{i}_0^{\vare}$ and 
$\mathbf{o}_0^{\vare}$ are divergence-free, the first two are tangent to $\Gamma_{\vare}$ 
and the last one vanishes on $\Gamma_{\vare}$. Moreover,
$\|\nabla \mathbf{o}_0^{\vare} \|_{L^2(\Pi_{\vare})}$ is bounded independent of $\vare$
and
\begin{equation} \label{gobi} 
\left\|\mathbf{o}_0^{\vare} - \alpha \varphi\left(\frac{\beta|x|}{\lambda_0}
\right) H\right\|_{L^2(\Pi_{\vare})} \to 
0 \mbox{ as } \vare \to 0,
\end{equation}
where $H = x^{\perp}/(2\pi|x|^2)$ and $\beta$ is as in Lemma \ref{confmap}. We also have that
\[\|\mathbf{i}_0^{\vare}\|_{L^2(\Pi_{\vare})} \leq C |\log \vare |^{\frac12}. \]
\end{lemma}  
  
\begin{proof}
Choose $\lambda_0$ such that the radii of the annulus where 
$\varphi^{\vare,\lambda_0}$ is not constant
are uniform in $\vare < \vare_0$. 

The $L^{\infty}$ bound on $\mathbf{b}^{\vare}_0$ comes from Theorem 4.1 in \cite{iln-shrink}.  
The $L^{p_1}$ bound on $\mathbf{b}^{\vare}_0$, $1 < p_1 < \infty$ follows from the local $L^{\infty}$ bound
above, together with two facts: (1) $m (1 - \varphi^{\vare,\lambda_0})
H^{\vare}$ has support in a compact set independent of $\vare$ and 
(2) $|K^{\vare}[\omega_0]| = \mathcal{O}(1/|x|^2)$ at infinity, uniformly in $\vare$.
As mentioned previously, fact (2) is  estimate (2.8) in \cite{iln-shrink}.

The $L^{p_3}$ bound on $\mathbf{o}^{\vare}_0$ follows from the fact that 
$\varphi^{\vare,\lambda_0}$ is constant outside an annulus independent of $\vare$, 
from formula \eqref{harm}, from the scaling $H^{\vare}(x)=\frac{1}{\vare} H_{\Pi}(x/\vare)$ and 
from the behavior of $T$ far from the obstacle given by Lemma \ref{confmap}.   

For $\mathbf{i}^{\vare}_0$, both the logarithmic estimate and the $L^{p_2}$ estimate follow
from adapting the argument used for estimate (3.7) of \cite{iln-shrink} in a straightforward
manner.

To estimate on $\nabla \mathbf{o}^{\vare}_0$ we observe that $|\nabla \mathbf{o}^{\vare}_0| = 
\mathcal{O}(1/|x|^2)$ near infinity, uniformly in $\vare$. This estimate easily reduces to 
an estimate  on $DH^{\vare}$, which in turn reduces to calculating derivatives of the 
conformal mapping $T$ using \eqref{holomform}.

Finally, \eqref{gobi} reduces to showing that $H-H^{\vare}$ goes to zero in $L^2$ near infinity,
which can be shown by a computation similar to the one carried out in the proof of Lemma 4.2 in \cite{iln-shrink}.   

\end{proof}  

In the remainder of this article, we will fix $\lambda_0$, independent of $\vare$, 
as in Lemma \ref{idp}, thereby fixing the bounded, inner and outer parts of 
the initial velocity, $\mathbf{b}_0^{\vare}$, $\mathbf{i}_0^{\vare}$ and $\mathbf{o}_0^{\vare}$.

Let us denote the Stokes semigroup on $\Pi_{\vare}$ by $S^{\vare}_{\nu}(t)$, 
so that $S^{\vare}_{\nu}(t)[v_0^{\vare}]$ is the solution to the Stokes 
system \eqref{stokes} on $\Pi_{\vare}$ with initial data $v_0^{\vare}$. We 
introduce the notation $\tau^{\vare} =  \tau^{\vare}(x) = \vare x$, the 
contraction by $\vare$. We observe the following fundamental relation between 
the Stokes system on $\Pi_{\vare}$ and on $\Pi$:
\begin{equation} \label{stokesrescale}
(S^1_{\nu}(t) [v_0^{\vare} \circ \tau^{\vare} ]) (x) = 
(S^{\vare}_{\nu}(\vare^2 t) [v_0^{\vare} ]) (\vare x), \mbox{ if } x \in \Pi.
\end{equation}

Our strategy to study the small obstacle limit begins by considering the solution $u^{\vare}$ 
of \eqref{NSevan} as a perturbation of $v^{\vare} \equiv S^{\vare}_{\nu}(t)u_0^{\vare}$.  The first 
thing we require is information on $v^{\vare}$ which we deduce in the result below.

\begin{lemma} \label{linest}
Let $\mathbf{b}^{\vare} \equiv S^{\vare}_{\nu}(t) \mathbf{b}_0^{\vare}$, $\mathbf{i}^{\vare} \equiv S^{\vare}_{\nu}(t) 
\mathbf{i}_0^{\vare}$, $\mathbf{o}^{\vare} \equiv S^{\vare}_{\nu}(t) \mathbf{o}_0^{\vare}$ and 
let $2 < p < \infty$. Then there exists a constant $K=K(p,\omega_0)>0$ such that for any $\vare>0$ we have:
\begin{enumerate}
\item[(i)]
$\displaystyle{\|\mathbf{b}^{\vare}\|_{L^p(\Pi_{\vare})} \leq K (\nu t)^{\frac{1}{p}-\frac{1}{2}},}$
\item[(ii)]
$\displaystyle{\|\mathbf{i}^{\vare}\|_{L^p(\Pi_{\vare})} \leq K|\gamma| (\nu t)^{\frac{1}{p}-\frac{1}{2}},}$
\item[(iii)]
$\displaystyle{\|\mathbf{o}^{\vare}\|_{L^p(\Pi_{\vare})} \leq K |\alpha| (\nu t)^{\frac{1}{p}-\frac{1}{2}}.}$
\end{enumerate}
\end{lemma} 

\begin{proof}
By \eqref{stokesrescale} we have that $\mathbf{b}^{\vare} (\vare^2 t, \vare x) = (S^1_{\nu} (t) [\mathbf{b}_0^{\vare} \circ \tau^{\vare}])(x)$, for $x \in \Pi$. Now, by Theorem \ref{ellpest}, item (S1), it follows that there exists $K_1>0$ such that  
\[\|S^1_{\nu} (t) [\mathbf{b}_0^{\vare} \circ \tau^{\vare}]  \|_{L^p(\Pi)} \leq K_1 (\nu t)^{\frac{1}{p}-\frac{1}{2}} \|\mathbf{b}_0^{\vare} \circ \tau^{\vare}\|_{L^2(\Pi)}.\]
Item (i) above follows from this estimate, together with \eqref{stokesrescale} and the fact that 
\[\|\mathbf{b}_0^{\vare} \circ \tau^{\vare} \|_{L^2(\Pi)} = \frac{1}{\vare} \|\mathbf{b}_0^{\vare}\|_{L^2(\Pi_{\vare})} \leq C\frac{1}{\vare},\]
where we have used Lemma \ref{idp} in the last inequality.
Items (ii) and (iii) follow in an analogous manner using Proposition \ref{stokesagapi} together with the fact that 
\[\|\mathbf{i}_0^{\vare} \circ \tau^{\vare} \|_{L^{2,\infty}(\Pi)} = \frac{1}{\vare} \|\mathbf{i}_0^{\vare}\|_{L^{2,\infty}(\Pi_{\vare})} \leq C\frac{|\gamma|}{\vare}\]
and \[\|\mathbf{o}_0^{\vare} \circ \tau^{\vare} \|_{L^{2,\infty}(\Pi)}  = \frac{1}{\vare} \|\mathbf{o}_0^{\vare}\|_{L^{2,\infty}(\Pi_{\vare})}\leq C\frac{|\alpha|}{\vare}.\]
We have used the scaling $H^{\vare} (x) = (1/\vare) H^1(x/\vare)$ above, see identity (3.6) in \cite{iln-shrink}.
\end{proof}

\vspace{0.5cm}

\begin{remark} \label{Kivare} {\em Using the rescaling \eqref{stokesrescale} we may deduce that the 
estimates (S1), (S2) and (S3) in Theorem \ref{ellpest} are valid in $\Pi_{\vare}$ with 
constants $K_1$, $K_2$ and $K_3$  independent of $\vare$.} 
\end{remark}

We will conclude this section with an observation on the amount of vorticity generated at the boundary in the initial
layer. This is a ``fixed $\vare$" calculation, before we take the vanishing obstacle limit. We denote 
vorticity associated to the velocity $u^{\vare}$, at time $t$, by $\omega^{\vare} = \omega^{\vare}(t,\cdot) \equiv 
\mbox{ curl }u^{\vare} (t,\cdot)$. Let us recall the discussion of flow in an exterior domain found in \cite{iln-shrink}. It was shown there that, if $\mbox{div }u^{\vare}(t,\cdot) = 0$, $\mbox{curl }u^{\vare}(t,\cdot) = \omega^{\vare}(t,\cdot)$, and if $u^{\vare}(t,\cdot)$ is tangent to $\Gamma_{\vare}$ and vanishes at infinity, then there exists unique $a=a(t)\in \real$ such that one can write $u^{\vare}(t,\cdot)$ as:
\[u^{\vare}(t,\cdot)=K^{\vare}[\omega^{\vare}(t,\cdot)] + a(t)H^{\vare};\]
(see Section 3.1 of \cite{iln-shrink} for details).
For the initial data \eqref{u0} we have, of course, $a(0)=\alpha$. We will prove that $a(t) = \alpha$ for any $t>0$. 
This fact relies on a result whose prove we defer to Section 8, see Corollary \ref{fixedvarevelest}. Using the notation introduced above we have:
\[ u^{\vare} - \mathbf{o}_0^{\vare} = K^{\vare}[\omega^{\vare}] + a(t)[1-\varphi^{\vare,\lambda_0} ]H^{\vare} + [a(t)  - \alpha] \varphi^{\vare,\lambda_0} H^{\vare}.\]
It will be proved in Corollary \ref{fixedvarevelest} that $u^{\vare} - \mathbf{o}_0^{\vare}$ belongs to 
$L^{\infty}_{\loc}([0,\infty);L^2(\Pi_{\vare}))$, although the estimate blows up as $\vare \to 0$. With this result in mind we first note that, from the estimates for $K^{\vare}$ proved in \cite{iln-shrink}, we expect to have that  
$K^{\vare}[\omega^{\vare}] \in L^{\infty}_{\loc}([0,\infty);L^2(\Pi_{\vare}))$. In fact, if $\omega^{\vare}$ were
compactly supported in space for each fixed time then this would follow from estimate (2.8) in \cite{iln-shrink};
the adaptation to our case is possible but escapes the scope of the present work. Next we recall that the harmonic vector field $H^{\vare}$ is smooth (because the conformal map $T^{\vare}$ is smooth and extends smoothly to $\Gamma^{\vare}$) and has $\mathcal{O}(1/|x|)$ behavior near infinity. Therefore we find that $[1-\varphi^{\vare,\lambda_0} ]H^{\vare} \in L^2(\Pi_{\vare})$, but $\varphi^{\vare,\lambda_0} H^{\vare} \notin L^2(\Pi_{\vare})$. Hence the only way for $u^{\vare} - \mathbf{o}_0^{\vare}$ to be square-integrable is for $a(t) - \alpha = 0$, as we wished. 

Since the flow $u^{\vare}$ satisfies the no-slip condition at any positive time, the circulation around
$\Gamma_{\vare}$ at $t>0$ vanishes. We make use once more of Stokes' Theorem to conclude that $0 = a(t) - m^{\vare}(t)$, where $m^{\vare}(t) = \int_{\Pi_{\vare}} \omega^{\vare}(t,x) \, dx$. We can now account precisely for the mass of vorticity produced at the boundary in the initial layer. We have: 
\[m^{\vare}(t) = \alpha.\]

\section{ Initial data asymptotics}

The purpose of this section is to study the limit, as $\vare \to 0$, of the initial velocity fields $u_0^{\vare}$. We begin by introducing some notation.

For each function $f$ defined on $\Pi_{\vare}$, we
introduce $Ef$, the extension of $f$ to $\real^2$, by setting $Ef \equiv 0$ in $\vare\Omega$.

\begin{lemma} \label{sobsilly} 
If $f \in W^{1,1}_{\loc}(\overline{\Pi}_{\vare})$ and if its trace vanishes on the 
boundary $\Gamma_{\vare}$ then $Ef \in W^{1,1}_{\loc}$ and $E\nabla f = \nabla Ef$.
\end{lemma}

The proof of this fact is elementary and we leave it to the reader. 

We will now introduce notation which will be used in the remainder of this paper.  
We denote by $\mathbb{P}$ the Leray projector on all of $\real^2$. Additionally, we introduce the cutoff 
\begin{equation} \label{etavare}
\eta^{\vare} = \eta^{\vare}(x) \equiv \varphi^{\vare,\vare}(x) = \varphi(|T^{\vare}(x)|)= \varphi(|T(x/\vare)|),
\end{equation}
where $T^{\vare}$, $\varphi^{\vare,\vare}$ and $\varphi$, were introduced in Section 4, 
see \eqref{tvare}, \eqref{cutoff}. Note that there exists a constant $C>0$ independent of $\vare$ 
such that $\eta^{\vare}(x) \equiv 1$ in $\{|x| > C\vare\}$.

Let 
\begin{equation} \label{biotker}
\mathcal{K} (x) = \frac{x^{\perp}}{2\pi|x|^2}
\end{equation}
be the kernel of the Biot-Savart law in all of $\real^2$, $f \mapsto \mathcal{K}\ast f$. Note that we 
denoted the same vector field by $H$ in \eqref{helena}. The different notations used for the same vector 
field are natural since $x^{\perp}/(2\pi|x|^2)$ plays two very different roles -- one as the 
kernel for the Biot-Savart law for the full plane and another as a harmonic generator for the cohomology of 
the exterior of any disk centered at the origin. 

Let $u_0^{\vare}$ be as in \eqref{u0} and $\gamma$ as in \eqref{mandgamma}.

\begin{lemma} \label{limidp}
Let $u_0 = \mathcal{K} \ast \omega_0 + \gamma H$. Then
we have that 
\[\mathbb{P}[\eta^{\vare} Eu_0^{\vare}] \to u_0 \mbox{ in } \mathcal{D}^{\prime}(\real^2)\] 
as $\vare \to 0$.
\end{lemma}

\begin{proof} 
We split $u_0^{\vare}$ in a different way than before:
\[u_0^{\vare} = (K^{\vare}[\omega_0] + m H^{\vare}) + \gamma H^{\vare} \equiv v_0^{\vare} + \gamma H^{\vare},\]
where $m$ was defined in \eqref{mandgamma}. By Lemma 4.2 in \cite{iln-shrink}, $\eta^{\vare}EH^{\vare} \to H$
strongly in $L^1_{\loc}(\real^2)$ as $\vare \to 0$, and by Lemma 4.1, in \cite{iln-shrink}  $\eta^{\vare}EH^{\vare}$
is divergence-free, so that $\mathbb{P} \eta^{\vare}EH^{\vare} = \eta^{\vare}EH^{\vare}$. Therefore,
\[\gamma\mathbb{P}[\eta^{\vare}EH^{\vare}] \to \gamma H \mbox{ in } \mathcal{D}^{\prime}(\real^2). \]

All that remains to prove is that 
$\mathbb{P}[\eta^{\vare} Ev_0^{\vare}] \to \mathcal{K} \ast \omega_0 \mbox{ in } \mathcal{D}^{\prime}(\real^2)$.
To see this, we begin by observing that $\eta^{\vare} Ev_0^{\vare}$ is uniformly bounded in $L^{\infty}(\real^2)$,
see Theorem 4.1 in \cite{iln-shrink}. Furthermore, we have additional control over the behavior of $Ev_0^{\vare}$
at infinity, so that there exists a constant $C>0$, independent of $\vare$, such
that $|\eta^{\vare} Ev_0^{\vare}| \leq C/|x|$. This follows from the explicit expressions for $K^{\vare}$,
$H^{\vare}$ given in (3.5) and (3.6) of \cite{iln-shrink}, from estimate (2.8) in \cite{iln-shrink}  
and from the compactness of the support of $\omega_0$. Therefore, $\eta^{\vare} Ev_0^{\vare}$ is
also uniformly bounded in $L^p(\real^2)$ for all $p>2$.  Fix $p > 2$
and let $\zeta \in L^{\infty}(\real^2) \cap L^p(\real^2)$ be a
weak-$\ast$ limit of $\{\eta^{\vare} Ev_0^{\vare}\}$.  

Next we observe that $\mbox{div }\eta^{\vare} Ev_0^{\vare} = \nabla \eta^{\vare} \cdot E v_0^{\vare}$
and $\mbox{curl }\eta^{\vare} Ev_0^{\vare} = \nabla^{\perp} \eta^{\vare} \cdot E v_0^{\vare} + 
\eta^{\vare} \omega_0$. The  cutoff $\eta^{\vare}$ is such that $|\nabla \eta^{\vare}|$ 
is bounded by $C/\vare$ and supported on a set of measure $C \vare^2$. Thus, $\nabla \eta^{\vare}$
(and $\nabla^{\perp} \eta^{\vare}$) converges to zero strongly in $L^q(\real^2)$ for any $1 \leq q<2$.
Hence, $\mbox{div }\eta^{\vare} Ev_0^{\vare} \to 0$ and 
$\mbox{curl }\eta^{\vare} Ev_0^{\vare} \to \omega_0$ strongly in $L^1$.  This, together with the
convergence of a subsequence of $\eta^{\vare} E v_0^{\vare}$ to $\zeta$ weak in $L^p$ implies that
$\mbox{div }\zeta = 0$ and $\mbox{curl }\zeta = \omega_0$ in the sense of distributions. 
Using that $\zeta \in L^p$ for $p<\infty$, we obtain that  $\zeta = \mathcal{K} \ast \omega_0$. 

Since we identified the limit, we have actually proved that 
$\eta^{\vare} E v_0^{\vare} \rightharpoonup \mathcal{K} \ast \omega_0$ weakly in $L^p$, without the need to
pass to subsequences. Therefore, as $\mathbb{P}$ is linear and continuous from $L^p$ to itself,
it follows that
\[\mathbb{P}[\eta^{\vare} E v_0^{\vare}] \rightharpoonup \mathbb{P}[\mathcal{K} \ast \omega_0] 
= \mathcal{K} \ast \omega_0,\]
which concludes the proof.

\end{proof}

\section{ The impulsively stopped rotating cylinder}

       In this section we will study the small obstacle asymptotics in the special case where
the obstacle is a small disk and the initial flow is harmonic. The purpose of this discussion is
just to motivate and illustrate our main result, as this special case is not rigorously required
for the remainder of the analysis.   

       Let us begin with a physical interpretation of our problem. Consider an infinite solid
cylinder of radius $r>0$ immersed in a viscous fluid occupying the whole space outside the 
cylinder. If the cylinder rotates with constant angular velocity $\lambda$, boundary friction
will induce rotational motion in the surrounding fluid which, one expects, will settle to a 
steady flow with velocity $u_0$ of the form  
\[u_0 = \frac{\lambda r^2 x^{\perp}}{|x|^2}, \]
see \cite{batchelor67} for a discussion of this example.

We consider viscous flow in the exterior of the cylinder with initial velocity $u_0$, imposing
the standard no-slip condition $u=0$ at $|x|=r$ for positive time. Physically, this corresponds 
to first ``preparing'' the initial data by rotating the cylinder for a long time, letting the
flow settle into the steady configuration $u_0$, and then suddenly halting the motion of the 
cylinder at time $t=0$. A shorthand description of this situation is that of the flow induced 
by an {\it impulsively stopped rotating cylinder}. The inconsistency between initial and 
boundary data in this problem generates a rather singular initial layer in the fluid motion. 

This problem has both translational symmetry, along the axis of the cylinder, and rotational 
symmetry about the same axis. If we assume that these symmetries are preserved by the flow then
the translational symmetry allows us to reduce the problem to two dimensions, and the rotational
symmetry cancels the nonlinearity and further reduces the dimension, so that the equations 
of motion reduce to the Stokes equation in the exterior of a disk. One may find 
in \cite{batchelor67} an explicit treatment of this problem, involving passing to polar coordinates, 
using separation of variables and expressing the solution by means of Fourier-Bessel integrals. 
        
The problem we wish to address next is the small obstacle limit of the 
impulsively stopped rotating cylinder as posed above. This means that
we consider $\Pi_{\vare} = \{|x| > \vare\}$. We have $H^{\vare} = x^{\perp}/(2\pi|x|^2)$, 
independent of $\vare$. In the notation of the previous section we pick $\omega_0 = 0$ and
\[u_0^{\vare} = \gamma \frac{x^{\perp}}{2\pi|x|^2}, \mbox{ in } \Pi_{\vare}.\]
Note that, $\gamma = 2 \pi \lambda \vare^2$, so that fixing the circulation $\gamma$ independent 
of $\vare$ means that the angular velocity $\lambda$ of the obstacle must blow up as the obstacle 
becomes smaller. 

We consider $u^{\vare} = u^{\vare}(x,t)$ and $p^{\vare} = p^{\vare}(x,t)$ solving 
\eqref{NSevan} with initial data as above. It is a nice exercise, which we leave to the reader, 
to prove that the solution preserves circular symmetry.  Preserving the symmetry means that
the velocity remains tangent to concentric circles about the origin, with the pressure and 
the magnitude of velocity both invariant under rotation.  One consequence of circular symmetry is
that $u^{\vare} \cdot \nabla u^{\vare}$ is a gradient field, so that we can absorb the nonlinearity
into the pressure term. Therefore, $u^{\vare}$ satisfies the Stokes system on $\Pi_{\vare}$.

We introduce the Lamb-Oseen vortex as the unique solution $U$ of 
\begin{equation} \label{linlimNS}
\left\{ \begin{array}{ll}
U_t + U\cdot \nabla U = -\nabla p + \nu \Delta U, & \mbox{ in } (0,\infty) \times \real^2 \\
\mbox{div }U = 0, & \mbox{ in } [0,\infty) \times \real^2 \\
\lim_{|x| \to \infty} U(t,x) = 0, & \mbox{ for all } t \geq 0\\
U(0,x) = \frac{x^{\perp}}{2\pi|x|^2}, 
\end{array} \right.
\end{equation}
see \cite{GGL,GW05}.
We are now ready to state and prove the main result in this section.

\begin{theorem} \label{linlim}
Let $\gamma \in \real$. Then the extension of velocity, $Eu^{\vare}$, converges weakly in 
$L^2_{\loc}((0,\infty),L^2_{\loc}(\real^2))$ to $\gamma U$.  
\end{theorem}  

We will present an outline of a proof for this result, highlighting the main ideas, for the 
sake of illustration. In Section 10 we will give a complete proof of our main theorem, which 
includes this example as a special case.  

\begin{proof}
By linearity we can assume, without loss of generality, that $\gamma = 1$. 

Let us begin by collecting the estimates for $Eu^{\vare}$. 
As in Lemma \ref{idp}, we decompose 
\[u^{\vare}_0 = \mathbf{b}_0^{\vare} + \mathbf{i}_0^{\vare} + \mathbf{o}_0^{\vare} \equiv 0 + (1 - \varphi) H^{\vare}
+ \varphi H^{\vare}, \]
with $\varphi$ the same cutoff used in \eqref{cutoff}.

We have:
\begin{enumerate}
\item $Eu^{\vare}$ is bounded in  $L^{q_0}_{\loc} ((0,\infty);L^p(\real^2))$ and in
$L^{q}_{\loc} ([0,\infty);L^p(\real^2))$, for any $2<p<\infty$, $q_0 = 2p/(p-2)$ and
$1\leq q <q_0$;
\item $Eu^{\vare}$ is bounded in $L^{\infty}_{\loc} ((0,\infty);L^2_{\loc}(\real^2))$;
\item $\nabla E u^{\vare}$ is bounded in $L^{2}_{\loc} ((0,\infty);L^2(\real^2))$
and in $L^{p}_{\loc} ([0,\infty);L^2(\real^2))$, for $ 1\leq p < 2$;
\item $Eu^{\vare} - \varphi H$ is bounded in $L^{\infty}_{\loc}((0,\infty);L^2(\real^2))$.
\end{enumerate}           

Indeed, we can use Lemma \ref{linest} to get the first item. 
The second and third items come from Lemma \ref{idp}, Theorem \ref{ellpest} and Proposition 
\ref{renen}. The last item comes from Lemma \ref{idp} and Proposition \ref{renen}. 

Let $\omega^{\vare} = \mbox{ curl }u^{\vare}$ be the vorticity. Let 
$\psi = \psi(x) \in C^{\infty}_c(\real^2)$  and $\theta = \theta(t) \in C^{\infty}_c((0,\infty))$
be test functions  and, for each $\delta > \vare$ consider $\varphi^{\delta} = \varphi(|x|/\delta)$, 
with $\varphi$ as before. The vorticity satisfies the vorticity equation 
in $\Pi_{\vare}$, which reduces to the heat equation in this special case, by symmetry.
Multiplying the vorticity equation by $\varphi^{\delta} \psi \theta$ and integrating by parts we
find
\begin{equation} \label{identitty}
\int_0^{\infty} \int_{\real^2} \varphi^{\delta} \psi \theta_t E\omega^{\vare} \, dxdt = - 
\nu  \int_0^{\infty} \int_{\real^2} \theta E\omega^{\vare} \Delta(\varphi^{\delta} \psi)\, dxdt. 
\end{equation}

Fix $E\omega^{\vare_k}$ a weakly converging subsequence in $L^2_{\loc}((0,\infty);L^2(\real^2))$
and let $\omega$ be its weak limit. Let $u$ be the weak limit of $Eu^{\vare_k}$, so that 
$\omega = \mbox{ curl }u$ for $t>0$. Note that all the a priori estimates (1-4) above pass to $u$ and $\omega$.
One may pass to the limit as $\vare_k \to 0$ in 
identity \eqref{identitty} to obtain

\begin{equation} \label{newtitty}
\int_0^{\infty} \int_{\real^2} \varphi^{\delta} \theta_t \psi \omega \, dxdt = -
\nu  \int_0^{\infty} \int_{\real^2} \theta\omega\Delta(\varphi^{\delta} \psi)\, dxdt. 
\end{equation}

Next, we assume that $\psi(0) = 0$.  

We have
\[\Delta (\varphi^{\delta} \psi) = \psi \Delta \varphi^{\delta} + 2 \nabla \psi \cdot \nabla \varphi^{\delta}
+ \varphi^{\delta} \Delta \psi. \]
Note that $\varphi^{\delta} -1 \to 0$ in $L^p$ for any $1 \leq p < \infty$, and that both $\nabla \varphi^{\delta}$
and $\psi \Delta \varphi^{\delta}$ converge weakly to zero in $L^2$ as
$\delta\to0$. The last convergence follows from the
hypothesis $\psi(0) = 0$. Therefore, taking the limit $\delta \to 0$ in \eqref{newtitty} we obtain

\begin{equation} \label{newertitty}
\int_0^{\infty} \int_{\real^2} \theta_t \psi \omega \, dxdt = -
\nu  \int_0^{\infty} \int_{\real^2} \theta\omega\Delta \psi \, dxdt. 
\end{equation}

This is enough to prove that the associated velocity $u$ is a weak solution of the Navier-Stokes
system in the full plane. To see that, take $\Phi$ a divergence-free, compactly supported test
vector field in the plane and define $\widetilde{\psi} \equiv \Delta^{-1} \mbox{ curl } \Phi$.
This means that $\Phi = \nabla^{\perp} \widetilde{\psi}$. For each $R>0$  let 
$\chi^R \equiv (1 - \varphi(|x|/R))$ and define
\[\Phi_R \equiv \nabla^{\perp}((\widetilde{\psi} - \widetilde{\psi}(0))\chi^R). \]
Take $\psi = \psi_R = (\widetilde{\psi} - \widetilde{\psi}(0))\chi^R$ in \eqref{newertitty}. We integrate
by parts to obtain
\begin{equation} \label{newestitty}
0 = \int_0^{\infty} \int_{\real^2} \theta_t \Phi_R u \, dxdt -
\nu  \int_0^{\infty} \int_{\real^2} \theta\omega\Delta \psi_R \, dxdt \equiv I_1 - \nu I_2. 
\end{equation}

We wish to pass to the limit $R\to \infty$ in both $I_1$ and $I_2$.
First note that
\[I_2 = \int_0^{\infty} \int_{\real^2} \theta\omega
[(\Delta \chi^R) (\widetilde{\psi} - \widetilde{\psi}(0)) + 2 \nabla\chi^R \cdot \nabla \widetilde{\psi}
+ \chi^R \Delta \widetilde{\psi} ] \, dxdt. \]
Observe that both $\Delta \chi^R$ and $\nabla \chi^R$ converge weakly to $0$ in $L^2$ and that $\chi^R \to 1$
pointwise. Therefore,
\[ \lim_{R\to\infty} I_2 = \int_0^{\infty} \int_{\real^2} \theta\omega\Delta \widetilde{\psi} \, dxdt. \]
On the other hand, since $\varphi H$ is independent of time, 
\[I_1 = \int_0^{\infty} \int_{\real^2} (\theta_t \Phi_R )(u - \varphi(|x|)H) \, dxdt \]
\[ = \int_0^{\infty} \int_{\real^2} \theta_t [(\nabla^{\perp}\chi^R)(\widetilde{\psi} - \widetilde{\psi}(0))
+ \chi^R \Phi] (u - \varphi(|x|)H) \, dxdt. \]
Therefore,
\[\lim_{R\to\infty} I_2 = \int_0^{\infty} \int_{\real^2} (\theta_t \Phi) (u-\varphi(|x|)H) \,dxdt
= \int_0^{\infty} \int_{\real^2} \theta_t \Phi u \,dxdt. \] 

Putting together the limits of $I_1$ and $I_2$, \eqref{newestitty} and integrating by parts 
in the vorticity term we obtain

\[\int_0^{\infty} \int_{\real^2} \theta_t \Phi u \, dxdt =
 - \nu  \int_0^{\infty} \int_{\real^2} \theta \Delta\Phi u \, dxdt.\]

This is the weak formulation of the velocity equation, under the
symmetry assumption. As for the initial data, it is not difficult to
prove that the initial data for the limit problem is the limit (in the
sense of distributions) of the initial data for the approximating
problems. This can be shown either by taking the test function $\theta$
not vanishing in 0 and including the initial data in the weak
fomulation and in the process of passing to the limit, or showing that
the velocity converges uniformly in time (up to time $t=0$) with values in
$H^{-3}_{\loc}$ in the same way as in Proposition \ref{tcomp}. 
Therefore, by Lemma \ref{limidp}, the velocity $u$ satisfies $u(0,x) =
H$. 
From uniqueness of weak solutions, 
see \cite{GGL,GW05}, the conclusion follows.      

\end{proof}

A result analogous to Theorem \ref{linlim} holds for solutions of the Stokes problem with arbitrary initial data, not necessarily symmetric, i.e., for any given $\gamma$, the small obstacle limit of solutions of the Stokes problem with initial data of the form \eqref{u0} is the solution of the Stokes problem in the full plane with initial data given by Lemma \ref{limidp}. The key issue here is linearity. As we will see in the next section, the smallness condition on the initial circulation is needed to deduce estimates in the initial layer, but it is solely due to the nonlinearity. 

Finally, we would like to use the impulsively stopped rotating cylinder as an illustration of what is taking place
in the initial layer. We describe the events at time $t=0$ for the flow generated by an impulsively stopped rotating cylinder noting that, for an (infinitesimally) small positive time, the fluid velocity vanishes at the boundary
but has a nonvanishing limit as one approaches the boundary from inside the fluid. Tangential discontinuities  
in fluid velocity are called {\it vortex sheets} in hydrodynamics. The effect of impulsively stopping the
rotation of the boundary amounts to placing a vortex sheet at the boundary and letting it diffuse into the 
bulk of the fluid through viscosity. The same rough picture describes what happens in the initial layer for 
our general problem as well.

\section{ Initial-layer and the nonlinear evolution}
    
We have fixed an arbitrary initial vorticity $\omega_0$ and circulation $\gamma$, and hence we must 
understand solutions of the Navier-Stokes equations with initial data which does not satisfy 
the no-slip boundary condition. The effect of the consequent initial layer can be understood roughly as
that of placing a vortex sheet at the boundary $\Gamma_{\vare}$ and letting it evolve, diffusing
into the flow. In Section 5 we saw that, if we consider linear evolution as described by the 
time-dependent Stokes semigroup, the a priori estimates from Section 2 are enough to establish
the small obstacle asymptotics. However, the problem of resolving this initial layer for the
full Navier-Stokes system and obtaining uniform estimates for the small obstacle problem is rather
more delicate and it is the subject of the present section.    

Let $u^{\vare}$ be the solution of the Navier-Stokes equations \eqref{NSevan} with initial
 velocity $u^{\vare}_0$ given by \eqref{u0} and let $v^{\vare} = S_{\nu}^{\vare}(t)[u^{\vare}_0]$ 
 as in Section 3. Let $W^{\vare} \equiv u^{\vare} - v^{\vare}$. Let $\mathbb{P}^{\vare}$ be the Leray 
projector on $\Pi_{\vare}$. The evolution of $W^{\vare}$ is described by the following system:
 \begin{equation} \label{wvareq}
 W^{\vare}_t - \nu \mathbb{P}^{\vare} \Delta W^{\vare} + \mathbb{P}^{\vare} \mbox{ div } (
 W^{\vare} \otimes W^{\vare}  + W^{\vare} \otimes v^{\vare} + v^{\vare} \otimes W^{\vare} +
 v^{\vare} \otimes v^{\vare}) = 0, 
 \end{equation}
 with the initial condition $W^{\vare}(0,x) = 0$ and the boundary condition $W^{\vare} = 0$ 
 on $\Gamma_{\vare}$. 
 
 We introduce the weighted-in-time norms. Let $p \geq 1$ and $f:(0,T) \to L^p(\Pi_{\vare})$ measurable. 
 Let $T>0$. We use the following notation:
 \[\|f\|_{p,T} \equiv \sup_{0\leq t \leq T} t^{\frac{1}{2} - \frac{1}{p}} \|f(t,\cdot)\|_{L^p(\Pi_{\vare})}.\]
 The use of these norms for the Navier-Stokes equations was pioneered
 by H. Fujita and T. Kato, see for example \cite{FK}.
 
 \begin{lemma} \label{initlayest} 
 Let $2 \leq p < \infty$. There exist positive constants $C_0$ and $C_p$  such that,
 if $0 < T \leq C_0 \nu^3$ and $|\gamma| < C_0 \nu$ then
 \[ \|W^{\vare}\|_{p,T} \leq C_p \nu^{\frac{p+2}{2p}}, \]
 for every $0 < \vare < \vare_0$.
 \end{lemma}   

 \begin{proof}
 We use Duhamel's principle to write
 \[W^{\vare}(t) = -\int_0^t S_{\nu}^{\vare}(t-\tau) \mathbb{P}^{\vare} \mbox{ div } (
 W^{\vare} \otimes W^{\vare}  + W^{\vare} \otimes v^{\vare} + v^{\vare} \otimes W^{\vare} +
 v^{\vare} \otimes v^{\vare})(\tau) \, d\tau.\]
 
 Take the $L^p$-norm and apply Theorem \ref{ellpest}, estimate (S3), with $2 \leq q \leq p$ to obtain
 \[\|W^{\vare}\|_{L^p(\Pi_{\vare})} \leq \int_0^t \| S_{\nu}^{\vare}(t-\tau) \mathbb{P}^{\vare} \mbox{ div } (
 W^{\vare} \otimes W^{\vare}  + W^{\vare} \otimes v^{\vare} + v^{\vare} \otimes W^{\vare} +
 v^{\vare} \otimes v^{\vare})(\tau)\|_{L^p(\Pi_{\vare})} \, d\tau \]
 \[ \leq K_3 \int_0^t (\nu (t-\tau))^{-\frac{1}{2} + \frac{1}{p} - \frac{1}{q}} \|
 (W^{\vare} \otimes W^{\vare}  + W^{\vare} \otimes v^{\vare} + v^{\vare} \otimes W^{\vare} +
 v^{\vare} \otimes v^{\vare})(\tau)\|_{L^q(\Pi_{\vare})} \, d \tau \]
 \[\leq K_3 \int_0^t (\nu (t-\tau))^{-\frac{1}{2} + \frac{1}{p} - \frac{1}{q}} 
 (\|W^{\vare}(\tau)\|_{L^{q_1}} \|W^{\vare}(\tau)\|_{L^{q_2}} + 
  \|W^{\vare}(\tau)\|_{L^{q_1}} \|v^{\vare}(\tau)\|_{L^{q_2}} + \] \[+  \|v^{\vare}(\tau)\|_{L^{q_1}}
 \|v^{\vare}(\tau)\|_{L^{q_2}}) \, d\tau, \]
where $q_1$ and $q_2$ are chosen so that $1/q = 1/q_1 + 1/q_2$ and we have used H\"{o}lder's inequality.
Next we use the definition of the $(p,t)$-norm to find
\[\|W^{\vare}\|_{L^p(\Pi_{\vare})} \leq K_3 ( \|W^{\vare}\|_{q_1,t} \|W^{\vare}\|_{q_2,t} + 
  \|W^{\vare}\|_{q_1,t} \|v^{\vare}\|_{q_2,t} + \] 
  \[+  \|v^{\vare}\|_{q_1,t} \|v^{\vare}\|_{q_2,t}) \int_0^t (\nu(t-\tau))^{
  -\frac{1}{2} + \frac{1}{p} - \frac{1}{q}} \tau^{\frac{1}{q}-1} \, d\tau. \]
  
We note that, for any $\alpha > -1$, $\beta > -1$, we have
\begin{equation} \label{emb1}
\int_0^t (t-\tau)^{\alpha} \tau^{\beta} \, d\tau \leq C(\alpha,\beta) t^{\alpha+\beta+1}.
\end{equation}
The proof of this inequality is an elementary calculation.

We wish to use \eqref{emb1} with $\alpha = -1/2 + 1/p - 1/q$ and $\beta = 1/q - 1$.
Assume that: 
\begin{equation} \label{phly}
2 \leq q \leq p.
\end{equation}
Note that this condition implies
\[\alpha = -\frac{1}{2} + \frac{1}{p} - \frac{1}{q} > -1 \hspace{1cm} 
\mbox{ and } \hspace{1cm} \beta = \frac{1}{q}-1 > -1.\]
Therefore, we find
\begin{equation} \label{nhoc} 
\|W^{\vare}\|_{p,t} \leq C \nu^{-\frac{1}{2} + \frac{1}{p} - \frac{1}{q}}
(\|W^{\vare}\|_{q_1,t} \|W^{\vare}\|_{q_2,t} + 
  \|W^{\vare}\|_{q_1,t} \|v^{\vare}\|_{q_2,t}  +  \|v^{\vare}\|_{q_1,t} \|v^{\vare}\|_{q_2,t})
\end{equation}
  
We divide the remainder of the proof in two steps: $p = 4$ and any $p \geq 2$.

First assume $p = 4$. Set $q_1 = q_2 = 4$, so that $q = 2$ and \eqref{phly} is satisfied. 
In this situation,
\eqref{nhoc} gives that $X(t) \equiv  \|W^{\vare}\|_{4,t}$ satisfies
\begin{equation} \label{poly}
X^2 + \left( \|v^{\vare}\|_{4,t} - \frac{\nu^{3/4}}{C} \right) X + \|v^{\vare}\|_{4,t}^2 \geq 0. 
\end{equation}

Note that $X(0) = 0$ and the parabola described by \eqref{poly} has, at $t=0$, two distinct nonnegative roots.
We observe that, {\it as long as this parabola has two distinct
  nonnegative roots $r_1(t) < r_2(t)$},  we have that inequality 
\eqref{poly} together with  the continuity in time of $X$, $r_1$ and $r_2$ imply that 
\begin{equation} \label{nhonhoc} 
0 \leq X(t) \leq r_1(t). 
\end{equation}  

The condition for the polynomial above to have two distinct roots is
\begin{equation} \label{discr0}
\left( \|v^{\vare}\|_{4,t} - \frac{\nu^{3/4}}{C} \right)^2 - 4 \|v^{\vare}\|_{4,t}^2 > 0.
\end{equation}

Since $\|v^{\vare}\|_{4,t} \geq 0$ and $\nu^{3/4}/C > 0$ we find that \eqref{discr0} is equivalent
to 
\begin{equation} \label{discr1}
\|v^{\vare}\|_{4,t} < \frac{\nu^{3/4}}{3C}. 
\end{equation}
Furthermore, under the above assumption, the two distinct roots are
also nonnegative.
The $(4,t)$-norm is nondecreasing in $t$ and hence, in order to guarantee that the polynomial in \eqref{poly}
have two distinct nonnegative roots, it is enough to verify \eqref{discr1} for $t=T$. Now we use the
linear estimates from the previous section to find conditions under which \eqref{discr1} is valid at $t=T$.

First, recall that $v^{\vare} = \mathbf{b}^{\vare} + \mathbf{o}^{\vare} + \mathbf{i}^{\vare}$, and that 
both $\mathbf{b}^{\vare}_0$ and $\mathbf{o}^{\vare}_0$ belong to $L^p(\Pi_{\vare})$, for $p>2$, with 
$L^p$-norms uniformly bounded in $\vare$, see Lemma \ref{idp}. We use estimate (S1) from
Theorem \ref{ellpest} with $p=q=4$, together with Lemma \ref{linest} to deduce:
\[\|v^{\vare}\|_{4,T} \leq \sup_{0 \leq t \leq T} t^{\frac{1}{4}} \|\mathbf{b}^{\vare}(t) +
\mathbf{o}^{\vare}(t)\|_{L^4(\Pi_{\vare})} + \sup_{0 \leq t \leq T} t^{\frac{1}{4}} \| 
\mathbf{i}^{\vare}(t)\|_{L^4(\Pi_{\vare})} \]
\[\leq K_1 T^{\frac{1}{4}}(\|\mathbf{b}_0^{\vare}\|_{L^4(\Pi_{\vare})}
+  \|\mathbf{o}_0^{\vare}\|_{L^4(\Pi_{\vare})}) + K|\gamma|\nu^{- \frac{1}{4}} 
\leq C (T^{\frac{1}{4}} +|\gamma|\nu^{- \frac{1}{4}}). \] 

Choose $C_0 >0$ so that the conditions 
\begin{equation} \label{conditions}
 T \leq C_0 \nu^{3} \hspace{1cm} \mbox{ and } \hspace{1cm} |\gamma| \leq C_0 \nu
\end{equation}
imply \eqref{discr1} with $t=T$.
 
 Assuming now that \eqref{conditions} are valid, using \eqref{nhonhoc} we have that
 \begin{equation} \label{estwvare}
  \|W^{\vare}\|_{4,t} = X(t) \leq r_1(t) \leq C \nu^{3/4}, 
 \end{equation}
 for $0 \leq t \leq T$. This concludes the proof in the case $p = 4$. 
 
 For any $p \geq 2$ we bootstrap the $(4,T)$-estimate in the following way.
 We return to \eqref{nhoc} and set $q_1 = q_2 = 4$.
 We then impose \eqref{conditions} to obtain
 
\[ \|W^{\vare}\|_{p,T} \leq C \nu^{-\frac{1}{2} + \frac{1}{p} -\frac{1}{2}}
(\|W^{\vare}\|_{4,T}^2+  \|W^{\vare}\|_{4,T} \|v^{\vare}\|_{4,T}  +  \|v^{\vare}\|_{4,T}^2)
\leq C(p)\nu^{\frac{p+2}{2p}}.\]
 
\end{proof}
 
 We conclude this section with the observation that
 $$\|W^{\vare}\|_{2,T} = \|W^{\vare}\|_{L^{\infty}((0,T);L^2(\Pi_{\vare}))},$$ so Lemma \ref{initlayest}
 actually provided a renormalized energy estimate on the initial layer.
 
\section{ Global-in-time nonlinear evolution}

In the previous section we obtained a priori estimates for $W^{\vare}$ in the initial layer
which are uniform in $\vare$. We will now splice the information we already possess with a 
standard energy estimate, in order to obtain a result which is global in time.
We retain the context introduced in the previous section.  

\begin{lemma} \label{globest}
Let $1 \leq p < 2$. Then $W^{\vare} \in   
L^{\infty}_{\loc}([0,\infty);L^2(\Pi_{\vare})) \cap L^{2}_{\loc}((0,\infty);H^1(\Pi_{\vare}))$,  
$W^{\vare} \in L^p_{\loc}([0,\infty);H^1(\Pi_{\vare}))$, and the respective norms are bounded independently of $\vare$.
\end{lemma}

\begin{remark} \label{notation} {\em Note that the bound in $L^{2}_{\loc}((0,\infty);H^1(\Pi_{\vare}))$ means that
$W^{\vare}$ is bounded in $L^{2}((\delta,T);H^1(\Pi_{\vare}))$ for any $0< \delta < T$, but not
necessarily for $\delta=0$.}
\end{remark}

\begin{proof}
We rewrite the evolution equation \eqref{wvareq} for $W^{\vare}$  as
\[ \left\{ \begin{array}{l}
W^{\vare}_t - \nu \Delta W^{\vare} +  (W^{\vare} + v^{\vare}) \cdot \nabla W^{\vare}  + W^{\vare} \cdot \nabla v^{\vare} + v^{\vare} \cdot \nabla v^{\vare} = -\nabla p, \mbox{ in } (0,\infty) \times \Pi_{\vare}
\\ \mbox{div }W^{\vare} = 0 \mbox{ in } [0,\infty) \times \Pi_{\vare} \\
W^{\vare}(0,\cdot) = 0 \mbox{ on } \{t=0\} \times \Pi_{\vare}\\ W^{\vare}(t,\cdot) = 0 \mbox{ on }
[0,\infty) \times \Gamma_{\vare}\end{array} \right. \]

We multiply the equation above by $W^{\vare}$ and integrate to obtain

\[ \mathcal{E} \equiv \frac{1}{2}\frac{d}{dt} \|W^{\vare}\|_{L^2}^2 + \nu \|\nabla W^{\vare}\|_{L^2}^2 = 
- \int_{\Pi_{\vare}} [W^{\vare}\cdot( W^{\vare} \cdot \nabla v^{\vare}) 
+ W^{\vare} \cdot (v^{\vare} \cdot \nabla v^{\vare})]\,dx \]
\[=  \int_{\Pi_{\vare}} [v^{\vare}\cdot( W^{\vare} \cdot \nabla W^{\vare}) 
+ v^{\vare} \cdot (v^{\vare} \cdot \nabla W^{\vare})]\,dx \]
\[\leq  \|W^{\vare}\|_{L^4} \|\nabla W^{\vare}\|_{L^2} \|v^{\vare}\|_{L^4} +
\|\nabla W^{\vare}\|_{L^2} \|v^{\vare}\|_{L^4}^2.\]

We will use the following interpolation inequality:
\[\|W^{\vare}\|_{L^4} \leq C \|W^{\vare}\|_{L^2}^{1/2} \|\nabla W^{\vare}\|_{L^2}^{1/2},\]
with a constant $C>0$ independent of $\vare$. This inequality in the case of $\real^2$  
can be found in Chapter 1 of \cite{lady69}. To obtain the corresponding inequality in
$\Pi_{\vare}$, one simply extends $W^{\vare}$ to $\real^2$ by setting it identically 
equal to zero inside $\vare\Omega$. As $W^{\vare}$ vanishes on $\Gamma_{\vare}$, the 
extension has $H^1$-norm in the plane identical to the $H^1$-norm of $W^{\vare}$ in
$\Pi_{\vare}$. Finally one uses the inequality in $\real^2$ on the extension.

We proceed with the estimate of $\mathcal{E}$:
\[\mathcal{E} \leq C \|W^{\vare}\|_{L^2}^{1/2} \|\nabla W^{\vare}\|_{L^2}^{3/2} \|v^{\vare}\|_{L^4} +
\|\nabla W^{\vare}\|_{L^2} \|v^{\vare}\|_{L^4}^2 \]
\[ \leq \frac{\nu}{2} \|\nabla W^{\vare}\|_{L^2}^2 + \frac{C}{\nu^3}\|W^{\vare}\|_{L^2}^2 \|v^{\vare}\|_{L^4}^4
+ \frac{1}{\nu} \|v^{\vare}\|_{L^4}^4, \]
where we used Young's inequality to estimate each of the products above. 
Next, we use Lemma \ref{linest} to deduce
\[ \|v^{\vare}\|_{L^4}^4 \leq \frac{C}{\nu t}. \]
Hence,
\[\frac{d}{dt} \|W^{\vare}\|_{L^2}^2 + \nu \|\nabla W^{\vare}\|_{L^2}^2 \leq
\frac{C}{\nu^4 t} \|W^{\vare}\|_{L^2}^2 + \frac{C}{\nu^2 t}, \]
for some constant $C$ independent of $\vare$.
Gronwall's inequality now gives, for any $0<t_1<t_2$,
\begin{equation} \label{enone}
\frac{\|W^{\vare}(t_2,\cdot)\|_{L^2}^2}{t_2^{C/\nu^4}} + \nu\int_{t_1}^{t_2} 
\frac{ \|\nabla W^{\vare}(s,\cdot)\|_{L^2}^2}{s^{C/\nu^4}} \, ds \leq \frac{\nu^2}{t_1^{C/\nu^4}} - \frac{\nu^2}{t_2^{C/\nu^4}} + \frac{\|W^{\vare}(t_1,\cdot)\|_{L^2}^2}{t_1^{C/\nu^4}}. 
\end{equation} 
First choose $t_1 = C_0 \nu^3/2$, with $C_0$ given in Lemma \ref{initlayest}. 
It follows from Lemma \ref{initlayest} with $p=2$ that
\[ \|W^{\vare}(t_1,\cdot)\|_{L^2}^2 \leq C \nu^2. \]
Therefore, 
\begin{equation} \label{gah1}
\|W^{\vare}(t,\cdot)\|_{L^2}^2 \leq C\nu^2 \left( \frac{Ct}{\nu^3} \right)^{C/\nu^4},
\end{equation}
for any $t \geq t_1$, and we conclude that $W^{\vare}$ is uniformly bounded in 
$L^{\infty}_{\loc}([0,\infty); L^2(\Pi_{\vare}))$ as desired.

Next, we return to \eqref{enone}, for the derivative estimate. Let $a>0$, multiply \eqref{enone}
by $t_1^{a+(C/\nu^4)-1}$ and integrate the resulting inequality with respect to $t_1$ from $0$ to $t_2$.
We obtain, 
\begin{equation} \label{gah2}
\int_0^{t_2} s^a \, \|\nabla W^{\vare}(s,\cdot)\|_{L^2}^2 \, ds \leq
\frac{a\nu^4 + C}{\nu^5} \left[ \frac{\nu^2}{a}t_2^a + \int_0^{t_2} 
s^{a-1} \, \|W^{\vare}(s,\cdot)\|_{L^2}^2 \, ds \right]. 
\end{equation}
Since we already know that $W^{\vare}$ is uniformly bounded in 
$L^{\infty}_{\loc}([0,\infty); L^2)$,
this estimate implies that $W^{\vare}$ is bounded in $L^2_{\loc}((0,\infty);H^1)$, uniformly in $\vare$.
Moreover, if $1\leq p < 2$, then the choice $a = (2-p)/2p$ above allows to conclude that  
$W^{\vare}$ is also bounded in $L^p_{\loc}([0,\infty);H^1)$.
\end{proof}

\section{ Velocity estimates}

In this section we derive global estimates on velocity using the analysis performed thus far. 
Before we begin, we require the following interpolation inequality. 

\begin{lemma} \label{interpgrr}
Let $2 < p < \infty$, $q_0 = 2p/(p-2)$, $1 \leq q \leq q_0$. Let $r \geq q(p-2)/p$
and set
\[\theta = \frac{2qr}{rp - q(p-2)}.\] 
If $p > 2q$ assume further that $r \leq q(p-2)/(p-2q)$. Then $\theta \geq 1$ and for any interval 
$I \subseteq \real$ and any $f \in L^r(I;H^1(\real^2)) \cap L^{\theta}(I;L^2(\real^2))$, 
we have
\[\|f\|_{L^q(I;L^p(\real^2))} \leq C \|f\|_{L^r(I;H^1(\real^2))}^{(p-2)/p} 
\|f\|_{L^{\theta}(I;L^2(\real^2))}^{2/p}. \]
\end{lemma}

\begin{proof}
We start by recalling the following standard interpolation inequality:
for any $g\in H^1(\real^2)$ we have
\begin{equation} \label{innerint}
\|g\|_{L^{2/(1-s)}(\real^2)} \leq C \|g\|_{H^s} \leq C \|g\|_{L^2}^{1-s}\|g\|_{H^1}^s, \mbox{ for any } 0 \leq s \leq 1.
\end{equation}

Fix exponents $p$, $q$, $r$ and $\theta$ as in the statement of this lemma. Observe that, if $ p > 2q$ then 
$\theta \geq 1$ if and only if $r \leq q(p-2)/(p-2q)$; we hence assume this further restriction on $r$ if $p > 2q$. 
In the other case, $p \leq 2q$, there is no additional restriction on $r$ to guarantee that $\theta \geq 1$.

Next, fix an interval $I \subseteq \real$ and let 
$f \in L^{r}(I;H^1(\real^2)) \cap L^{\theta}(I;L^2(\real^2))$. Let $s = (p-2)/p$, so that $2/(1-s) = p$. We use \eqref{innerint} and H\"{o}lder's inequality to obtain:
\begin{eqnarray*}
\|f\|_{L^q(I;L^{p}(\real^2))}^q = \int_I \|f(\tau,\cdot)\|_{L^p}^q\, d\tau \leq C \int_I \|f(\tau,\cdot)\|_{H^1}^{q(p-2)/p} \|f(\tau,\cdot)\|_{L^2}^{2q/p} \, d\tau 
\\
\leq C \|f\|_{L^{r}(I;H^1(\real^2))}^{q(p-2)/p} \|f\|_{L^{\theta}(I;L^2(\real^2))}^{2q/p},
\end{eqnarray*}
which concludes the proof. The condition $r \geq q(p-2)/p$ was used in 
H\"{o}lder's inequality when estimating the product of two functions in $L^{rp/q(p-2)}(I)$ and $L^{\theta p/(2q)}(I)$ 
above, so as to guarantee that  $rp/q(p-2) \geq 1$. 

\end{proof}

\begin{theorem} \label{velest}
Let $u^{\vare}$ be the solution of \eqref{NSevan} with initial velocity 
$u_0^{\vare}$ as in \eqref{u0} and recall that $u_0^{\vare} = \mathbf{b}_0^{\vare} + 
\mathbf{i}_0^{\vare} + \mathbf{o}_0^{\vare}$. Then the following hold true. 
\begin{enumerate}
\item Let $2<p<\infty$, $q_0 = 2p/(p-2)$ and $1\leq q < q_0$. Then $\{Eu^{\vare}\}$ is bounded 
in $L^{q_0}_{\loc}((0,\infty);L^p) \cap L^q_{\loc}([0,\infty);L^p)$.
\item The family $\{Eu^{\vare} - E\mathbf{o}_0^{\vare}\}$ is bounded in 
$L^{\infty}_{\loc}((0,\infty);L^2)$ and the family $\{Eu^{\vare} - E\mathbf{o}_0^{\vare} - E\mathbf{i}^{\vare}\}$
is bounded in $L^{\infty}_{\loc}([0,\infty);L^2)$.
\item The family $\{Eu^{\vare}\}$ is bounded in  $L^{\infty}_{\loc}((0,\infty);L^2_{\loc})$.
\item For any $1 \leq p < 2$, we have $\{\nabla E u^{\vare}\}$ is bounded in 
$L^2_{\loc}((0,\infty);L^2) \cap L^p_{\loc}([0,\infty);L^2)$.
\end{enumerate}
\end{theorem} 

\begin{proof}
Statement (1) involves two estimates: the first one on the open time interval $(0,\infty)$ and the
second on the closed interval $[0,\infty)$. We begin by addressing the first estimate.  

Fix $2 < p < \infty$. Fix $0< \delta <T$ and set $I = (\delta,T)$.  We first show that 
$\{Eu^{\vare} - E\mathbf{o}_0^{\vare}\}$ is bounded in  
$L^{\infty}((\delta,T);L^2(\real^2)) \cap L^2((\delta,T);H^1(\real^2))$.
We write
\begin{equation} \label{decompp}
Eu^{\vare} - E\mathbf{o}_0^{\vare} = E(\mathbf{i}^{\vare} + \mathbf{b}^{\vare}) + E(\mathbf{o}^{\vare} - \mathbf{o}_0^{\vare}) +
EW^{\vare} \equiv A_1 + A_2 + A_3. 
\end{equation}

We observe that $A_1$ is bounded in $L^{\infty}(I;H^1(\real^2))$. To see that, choose $1 < r <2$
and use Theorem \ref{ellpest} together with Lemma \ref{idp} and Remark \ref{Kivare}   to obtain
\begin{equation} \label{duh}
\begin{array}{r}
t^{\frac{1}{r}-\frac{1}{2}} \left( \|\mathbf{i}^{\vare}\|_{L^2(\Pi_{\vare})} + \|\mathbf{b}^{\vare}\|_{L^2(\Pi_{\vare})} \right) + t^{\frac{1}{r}} \left( \|\nabla \mathbf{i}^{\vare}\|_{L^2(\Pi_{\vare})} + \| \nabla \mathbf{b}^{\vare}\|_{L^2(\Pi_{\vare})} 
\right) \\ 
\leq (K_1 + K_2) \left( \|\mathbf{i}^{\vare}_0\|_{L^r(\Pi_{\vare})} + \|\mathbf{b}^{\vare}_0\|_{L^r(\Pi_{\vare})} 
\right) \leq K(r), 
\end{array} 
\end{equation}
for some $K(r) > 0$, independent of $\vare$. The estimate on $A_1$ follows from Lemma \ref{sobsilly}, together
with the inequality above. For $A_2$, we use Proposition \ref{renen}, together with Lemma \ref{idp} to
conclude that 
\begin{equation} \label{duh2}
\|\mathbf{o}^{\vare}(t,\cdot) - \mathbf{o}_0^{\vare}\|_{L^2(\Pi_{\vare})}^2 + 
\nu \int_0^t  \|\nabla \mathbf{o}^{\vare}(s,\cdot) - \nabla \mathbf{o}_0^{\vare}\|_{L^2(\Pi_{\vare})}^2 \, ds
\leq \nu t \|\nabla \mathbf{o}_0^{\vare}\|_{L^2(\Pi_{\vare})}^2 \leq C, 
\end{equation}
for some $C>0$ independent of $\vare$. This, together with Lemma \ref{sobsilly}, implies that $A_2$ 
is uniformly bounded in $L^{\infty}(I;L^2(\real^2)) \cap L^2(I;H^1(\real^2))$. 
For the estimate on $A_3$, we simply use Lemma \ref{globest}, together with Lemma \ref{sobsilly},
showing that $A_3$  is uniformly bounded in $L^{\infty}(I;L^2(\real^2)) \cap L^2(I;H^1(\real^2))$ as well.

We use Lemma \ref{interpgrr} with $q = q_0$ and $r=2$, so that $\theta = \infty$, to conclude that
$\{Eu^{\vare} - E\mathbf{o}_0^{\vare}\}$ is bounded in $L^{q_0}(I;L^p(\real^2))$.
Next we note that $\{ E\mathbf{o}_0^{\vare} \}$ is uniformly bounded in $L^p(\real^2)$
by Lemma \ref{idp}, which concludes this portion of the proof.

We now address the second part of (1), which is an estimate on the closed time interval $[0,\infty)$.
The difficulty here is that we do not have Leray-type estimates on the pieces of $u^{\vare}$ all the way
down to $t=0$, so that the result becomes more delicate.

Fix $2 < p < \infty$ and $1 \leq q < q_0$. Let $T>0$ and set $I = [0,T]$.
We consider again the decomposition \eqref{decompp} and we estimate each piece. 
Estimate \eqref{duh}, together with Lemma \ref{sobsilly} implies that $A_1$ is uniformly
bounded on $L^{r_1}(I;H^1(\real^2)) \cap L^{\theta_1}(I;L^2(\real^2))$, for any $1 \leq r_1 < 2$
and any $1 \leq \theta_1 < \infty$. We use Lemma \ref{interpgrr} with $p$ and $q$ as above. We need
to find $r \in [1,2)$ satisfying the restrictions in Lemma \ref{interpgrr} in order to be able 
to use $r=r_1$. This is always possible because the restriction on $r$ always includes $r \geq q(p-2)/p$,
and $q < q_0$ is equivalent to $q(p-2)/p < 2$. This implies that $A_1$ is bounded in $L^q(I;L^p(\real^2))$.
For $A_2$, we merely observe that \eqref{duh2} gives an uniform bound in $L^{\infty}(I;L^2(\real^2)) 
\cap L^2(I;H^1(\real^2))$ on $A_2$, which in turn yields the desired estimate. To treat $A_3$, we put together
Lemma \ref{globest} and Lemma \ref{sobsilly} to conclude that $A_3$ is uniformly bounded in 
$L^{\infty}(I;L^2(\real^2)) \cap L^{r_2}(I;H^1(\real^2))$, for any $1 \leq r_2 < 2$. Clearly, this is enough
to obtain  the estimate in $L^q(I;L^p(\real^2))$ for $A_3$. The proof of (1) is concluded once we recall
the observation that $\{ E\mathbf{o}_0^{\vare} \}$ is uniformly bounded in $L^p(\real^2)$, which we already 
used in the proof of the first part of (1). 

We now address statement (2), which also consists of two estimates. The proof of the first estimate in (2) is 
contained in the proof of the first part of (1), the estimate on the open time interval. As for the second 
estimate in item (2), we write 
\[Eu^{\vare} - E\mathbf{o}_0^{\vare} - E\mathbf{i}^{\vare} = E\mathbf{b}^{\vare} + E(\mathbf{o}^{\vare}-\mathbf{o}_0^{\vare})
+ EW^{\vare}.\]
We have already shown that the second and third terms in the decomposition above are bounded in 
$L^{\infty}_{\loc}([0,\infty);L^2(\real^2))$. The first term satisfies 
\[\|E\mathbf{b}^{\vare}(t,\cdot)\|_{L^2} \leq  \|\mathbf{b}^{\vare}(t,\cdot)\|_{L^2(\Pi_{\vare})} \leq 
C\|\mathbf{b}^{\vare}_0\|_{L^2(\Pi_{\vare})} \leq C, \]
by Theorem \ref{ellpest} and Lemma \ref{idp}. 

The third item, statement (3), can be obtained from (2) by observing that, by Lemma \ref{idp} and 
Lemma \ref{sobsilly},  $E\mathbf{o}_0^{\vare}$ is uniformly bounded in $L^r(\real^2)$, for any $r>2$, 
which is contained in $L^2_{\loc}(\real^2)$. 

Statement (4) again consists of two estimates, one on the open time interval, the other on the closed interval. 
The estimates on the open time interval are trivially contained in the proof of the
first estimate in item (1), once we observe that Lemma \ref{idp} and Lemma \ref{sobsilly} give a uniform estimate
in $L^2(\real^2)$ for $\nabla E \mathbf{o}_0^{\vare}$. Similarly, the proof of the second part of (4) is contained
in the proof of the second part of (1), together with the $L^2(\real^2)$ estimate for $\nabla E \mathbf{o}_0^{\vare}$
which we have just derived. 

This concludes the proof.
 
\end{proof}

The last result in this section is an estimate for ``fixed $\vare$", which was already used to deduce that the amount of vorticity generated at the boundary in the initial layer, $m^{\vare}(t)=\int_{\Pi_{\vare}} \omega^{\vare}(t,x)\, dx$ is the same for any $\vare$ and equals $\alpha$.

\begin{corollary} \label{fixedvarevelest}
For each fixed $0 < \vare < \vare_0$ we have: $Eu^{\vare} - E\mathbf{o}_0^{\vare} \in L^{\infty}_{\loc}([0,\infty);L^2(\Pi_{\vare}))$. 
\end{corollary}

\begin{proof}
This is an immediate consequence of item (2) in Theorem \ref{velest}, together with the observation that $E\mathbf{i}^{\vare} \in L^{\infty}_{\loc}([0,\infty);L^2(\Pi_{\vare}))$ for each fixed $\vare$. This last fact follows from the estimates on the Stokes semigroup in Theorem \ref{ellpest} and the fact that $E\mathbf{i}_0^{\vare} \in L^2(\Pi_{\vare})$ with an $L^2$-norm that blows up as $|\log \vare|^{\frac12}$, see Lemma \ref{idp}.
\end{proof}

\section{ Compactness in space-time}

As far as a priori estimates go, the last ingredient we require is uniform control on how solutions 
evolve in time. This is often very easy to accomplish once spatial estimates are in place because
ultimately the PDE itself is nothing more than an expression of time-derivatives of the solution in
terms of spatial information. In our case, however, there are two difficulties that will make
this step of the analysis somewhat involved: (i) the spatial estimates available are for $Eu^{\vare}$,
which does not satisfy a PDE: and (ii) the nonlinearity in our problem is quadratic, but the estimate
up to time zero in Theorem \ref{velest} item (1) is $L^p$, $p>2$, which entails problems at infinity.  
We deal with these difficulties through the following main ideas: we use the vorticity equation to 
describe the time evolution, we use the interplay of vorticity and velocity and we renormalize problem 
terms. 

  Let $\Phi$ be a smooth, compactly supported vector field and consider $\mathbb{P}$ the Leray
projector for the plane. We consider the Hodge decomposition of the vector field $\Phi$, given by
$\Phi = \mathbb{P} \Phi + (\Phi - \mathbb{P} \Phi)$. The divergence-free part $\mathbb{P}\Phi$ 
is smooth, but not compactly supported. In fact,  
\begin{equation} \label{decay}
|\mathbb{P}\Phi(x)| = \mathcal{O}\left(\frac{1}{|x|^2}\right), \mbox{ as } |x| \to \infty,
\end{equation}
see the proof of Proposition 1.16 in \cite{MB02}. Let $\psi = \psi(x)$ be the stream function
associated with $\mathbb{P} \Phi$, so that $\nabla^{\perp}\psi = \mathbb{P} \Phi$. We assume
that $\psi(0)=0$, at the expense of having $\psi = \mathcal{O}(1)$ at $\infty$. Clearly,
$|\psi(x)| \leq |x| \|\nabla \psi\|_{L^{\infty}}$, so that using the Sobolev imbedding 
$H^2 \hookrightarrow L^{\infty}$ followed by the fact that 
$\nabla \psi = \nabla (\Delta)^{-1} \nabla^{\perp} \cdot \Phi$, a zeroth order singular integral 
operator acting on $\Phi$, we obtain
\begin{equation} \label{veryclever}
|\psi(x)| \leq C |x| \|\Phi\|_{H^2(\real^2)}.
\end{equation}
 
We now observe that for $1<q<\infty$ and  $\chi \in
C^{\infty}_c(\real^2)$, there exists a
constant $C = C(q,\chi) > 0$ such that 
\begin{equation} \label{Biotest}
\|\mathbb{P}[\chi \Phi]\|_{L^q(\real^2)} \leq C
\|\Phi\|_{H^2(\real^2)}. 
\end{equation}
Indeed, as $\mathbb{P}$ is a zero-th order singular integral operator we have
\[\|\mathbb{P}[\chi \Phi]\|_{L^q(\real^2)} \leq C\|\chi \Phi\|_{L^q(\real^2)}
\leq C\|\chi\|_{L^q(\real^2)} \|\Phi\|_{L^\infty(\real^2)} \leq C(q,\chi)\|\Phi\|_{H^2(\real^2)}.\]

Recall the extension operator $E$, introduced in beginning of Section 4.
For each $\vare>0$, consider the cutoff $\eta^{\vare}$ introduced in \eqref{etavare}.
We are ready to state and prove the main result in this section. 

\begin{proposition} \label{tcomp}
The sequence $\{\mathbb{P}[\eta^{\vare}Eu^{\vare}]\}$ is precompact in $L^{\infty}_{\loc}([0,\infty);H^{-3}_{\loc}(\real^2))$.
\end{proposition}    

\begin{proof}
Fix $\Phi$ a smooth, compactly supported vector field and let 
\[\psi = \psi(x) = [(\Delta)^{-1} \mbox{ curl } \Phi](x) -  [(\Delta)^{-1} \mbox{ curl } \Phi](0),\]
satisfying \eqref{veryclever}. 
For each $t \geq 0$, we introduce  an auxiliary functional $F^{\vare} = F^{\vare}(t) \in H^{-2}(\real^2)$ 
defined by
\[\langle F^{\vare}(t),\Phi \rangle = \int (Eu^{\vare}(t,x) - E\mathbf{o}_0^{\vare}(x)) \cdot 
(\nabla^{\perp}\eta^{\vare})(x) \psi(x) \,dx. \]
The proof will be divided into two steps. We will show that, for each $(t_1,t_2) \subset [0,\infty)$, we have  
$\{\mathbb{P}[\eta^{\vare}Eu^{\vare}]+ F^{\vare}\}$ is bounded and equicontinuous as
a function of $(t_1,t_2)$ into $H^{-2}_{\loc}(\real^2)$ and we will show that $F^{\vare} \to 0$ strongly in
$L^{\infty}_{\loc}([0,\infty);H^{-2}(\real^2))$. The desired conclusion follows from these two steps by 
using Arzela-Ascoli's Theorem. 

Let us begin by proving that $F^{\vare} \to 0$ strongly in $L^{\infty}_{\loc}([0,\infty);H^{-2}(\real^2))$. Indeed,
we use Theorem \ref{velest}, \eqref{veryclever} and the properties of the cutoff $\eta^{\vare}$ to deduce
\begin{equation} \label{fvare}
 |\langle F^{\vare}(t),\Phi \rangle| \leq \|Eu^{\vare}(t,\cdot) - E\mathbf{o}_0^{\vare}\|_{L^2} 
\|\psi \nabla^{\perp} \eta^{\vare}\|_{L^2} 
\end{equation}
\[ \leq C \|Eu^{\vare}(t,\cdot) - E\mathbf{o}_0^{\vare}\|_{L^2} \left(\int_{|x|<C\vare} |x|^2 \|\Phi\|_{H^2}^2 |\nabla \eta^{\vare}|^2 \, dx \right)^{1/2} \] 
\[ \leq C \vare\|Eu^{\vare}(t,\cdot) - E\mathbf{o}_0^{\vare}\|_{L^2} \|\Phi\|_{H^2} \]
\[\leq C \vare (\|E\mathbf{i}^{\vare}(t,\cdot)\|_{L^2} + \|Eu^{\vare}(t,\cdot) - E\mathbf{o}_0^{\vare} -
E\mathbf{i}^{\vare}(t,\cdot)\|_{L^2} ) \|\Phi\|_{H^2} 
\leq C \vare (C_1|\log\vare|^{\frac12}+ C_2) \|\Phi\|_{H^2},\]
by Lemma \ref{idp}, Theorem \ref{ellpest} and Theorem \ref{velest}. Clearly, this proves our assertion. 
The proof of the other assertion is a bit more involved.

We introduce a cutoff for infinity. For each $R>0$, let $\chi^R = \chi^R(x) = 1 - \varphi(|x|/R)$.   
The vorticity $\omega^{\vare} = \mbox{ curl }u^{\vare}$ satisfies the equation
\[\omega^{\vare}_t + u^{\vare} \cdot \nabla \omega^{\vare} = \nu \Delta \omega^{\vare}. \]
Let $0 \leq t_1 < t_2 < \infty$ and denote the interval $[t_1,t_2]$ by $J$. 
We multiply the vorticity equation by $\eta^{\vare} \psi \chi^R$, 
integrate in space and time between $t_1$ and $t_2$, and integrate by parts to obtain
\begin{equation} \label{lhs}
\int [Eu^{\vare}(t_2,\cdot) - Eu^{\vare}(t_1,\cdot)] \cdot \nabla^{\perp}(\eta^{\vare} \psi \chi^R)\, dx 
\end{equation}
\[ = \int_{t_1}^{t_2} \int (Eu^{\vare} \cdot \nabla E\omega^{\vare}) \eta^{\vare} \psi \chi^R \, dxdt
- \nu \int_{t_1}^{t_2} \int (\Delta E\omega^{\vare}) \eta^{\vare} \psi \chi^R \, dxdt \equiv I_1 - I_2. \]
We first estimate $I_1$. We integrate by parts and deduce
\[ I_1 = -\int_{t_1}^{t_2} \int (Eu^{\vare} \cdot \nabla\eta^{\vare}) 
E\omega^{\vare} \psi \chi^R \, dxdt + \int_{t_1}^{t_2} \int (Eu^{\vare} \cdot (\mathbb{P}\Phi)^{\perp}) 
E\omega^{\vare} \eta^{\vare} \chi^R \, dxdt\] 
\[- \int_{t_1}^{t_2} \int (Eu^{\vare} \cdot \nabla\chi^R )
E\omega^{\vare} \psi \eta^{\vare} \, dxdt \equiv -I_{11} + I_{12} - I_{13}.\]
Using H\"{o}lder's inequality first in space and then in time we have
\[|I_{11}| \leq \int_{t_1}^{t_2} \|Eu^{\vare}\|_{L^4} \|E\omega^{\vare}\|_{L^2} \|\psi \nabla \eta^{\vare}\|_{L^4}
\|\chi^R\|_{L^{\infty}} \, dt \]
\[\leq C \|\psi \nabla \eta^{\vare}\|_{L^4} \|Eu^{\vare}\|_{L^3(J;L^4)} 
\|E\omega^{\vare}\|_{L^{9/5}(J;L^2)} |t_2-t_1|^{1/9} \leq C \|\Phi\|_{H^2} |t_2 - t_1|^{1/9}, \]
where in the last inequality we used Theorem \ref{velest} and we used again \eqref{veryclever} along 
with the properties of $\eta^{\vare}$. 

Similarly, we have
\[|I_{12}| \leq \|\mathbb{P} \Phi \|_{L^4} \|Eu^{\vare}\|_{L^3(J;L^4)} 
\|E\omega^{\vare}\|_{L^{9/5}(J;L^2)} |t_2 -t_1|^{1/9} \leq C \|\Phi\|_{H^2} |t_2 - t_1|^{1/9}, \]
since $\mathbb{P}$ is a zeroth order operator and $H^2 \hookrightarrow L^4$. 
Finally,
\[|I_{13}| \leq C \|\psi \eta^{\vare}\|_{L^{\infty}} \|\nabla
\chi^R\|_{L^4} 
\|Eu^{\vare}\|_{L^3(J;L^4)} 
\|E\omega^{\vare}\|_{L^{2/3}(J;L^2)}
\leq C \|\Phi\|_{H^2}R^{-1/2}, \]
as 
\[|\psi(\infty)| = |[(\Delta)^{-1} \mbox{ curl }\Phi](0)| \leq C\|\Phi\|_{H^2}\] 
and $\nabla \chi^R = \mathcal{O}(1/R)$, supported on a set of measure $\mathcal{O}(R^2)$. 

Therefore, 
\begin{equation} \label{bag1}
\limsup_{R\to \infty} |I_1| \leq C \|\Phi\|_{H^2} |t_2-t_1|^{1/9}.
\end{equation}

Next we treat $I_2$. We integrate by parts and use the fact that the support of $\nabla \eta^{\vare}$ and
of $\nabla \chi^R$ are disjoint to obtain:
\[I_2 = \nu \int_{t_1}^{t_2} \int E\omega^{\vare} (\Delta \eta^{\vare})\psi \chi^R \, dxdt +
\nu \int_{t_1}^{t_2} \int E\omega^{\vare} \eta^{\vare} (\Delta \psi) \chi^R \, dx dt\]
\[+ \nu \int_{t_1}^{t_2} \int E\omega^{\vare} \eta^{\vare} \psi (\Delta \chi^R) \, dxdt 
- 2 \nu \int_{t_1}^{t_2} \int E\omega^{\vare} (\chi^R \nabla \eta^{\vare} + \eta^{\vare}\nabla\chi^R) 
(\mathbb{P} \Phi)^{\perp} \, dxdt\] 
\[= \nu I_{21} + \nu I_{22} + \nu I_{23} - 2 \nu I_{24}. \]
By arguments similar to those used for $I_1$ we have
\[|I_{21}| \leq \|\Delta \eta^{\vare}\psi\|_{L^2} \|E\omega^{\vare}\|_{L^{9/5}(J;L^2)} |t_2 - t_1|^{4/9} \leq
C \|\Phi\|_{H^2}|t_2 - t_1|^{4/9};\]
\[|I_{22}| \leq \|\Delta \psi\|_{L^2} \int_{t^1}^{t_2} \|E\omega^{\vare}\|_{L^2} \, dt \leq 
C \|\mbox{ curl } \Phi\|_{L^2} |t_2 - t_1|^{4/9} \leq C \|\Phi\|_{H^2}|t_2 - t_1|^{4/9};\]
\[|I_{23}| \leq \|\Delta \chi^R\|_{L^2} \|\psi\|_{L^{\infty}} \int_{t_1}^{t_2} \|E\omega^{\vare}\|_{L^2} \,dt
\leq C \|\Phi\|_{H^2};\] 
\[|I_{24}| \leq \|\mathbb{P} \Phi\|_{L^{\infty}} \|\chi^R \nabla\eta^{\vare} + \eta^{\vare} \nabla \chi^R\|_{L^2}
\int_{t_1}^{t_2} \|E\omega^{\vare}\|_{L^2} \, dt \leq C \|\Phi\|_{H^2} |t_2 - t_1|^{4/9}. \] 

Therefore,
\begin{equation} \label{bag2}
\limsup_{R\to \infty} |I_2| \leq C \|\Phi\|_{H^2} |t_2-t_1|^{4/9}.
\end{equation}

 We expand the left hand side of identity \eqref{lhs} to find

\[\int [Eu^{\vare}(t_2,\cdot) - Eu^{\vare}(t_1,\cdot)] \cdot \nabla^{\perp}(\eta^{\vare} \psi \chi^R)\, dx \]
\[= \int [Eu^{\vare}(t_2,\cdot) - Eu^{\vare}(t_1,\cdot)] \cdot (\nabla^{\perp}\eta^{\vare}) \psi \chi^R\, dx
+ \int [Eu^{\vare}(t_2,\cdot) - Eu^{\vare}(t_1,\cdot)] \cdot \eta^{\vare}(\nabla^{\perp}\psi) \chi^R\, dx\]
\[+ \int [Eu^{\vare}(t_2,\cdot) - Eu^{\vare}(t_1,\cdot)] \cdot \eta^{\vare} \psi (\nabla^{\perp} \chi^R) \, dx
= A_1 + A_2 + A_3. \]  

We will show that each of the $A_i$'s  has a limit when $R \to \infty$. To see that,
first note that 
\[Eu^{\vare}(t_2,\cdot) - Eu^{\vare}(t_1,\cdot) \] 
\[ =  [Eu^{\vare}(t_2,\cdot) - E\mathbf{o}_0^{\vare} - E\mathbf{i}^{\vare}(t_2,\cdot)] - [Eu^{\vare}(t_1,\cdot) - E\mathbf{o}_0^{\vare}- E\mathbf{i}^{\vare}(t_1,\cdot)] + E\mathbf{i}^{\vare}(t_2,\cdot) - E\mathbf{i}^{\vare}(t_1,\cdot), \]
which belongs to $L^2(\real^2)$, for each fixed $\vare > 0$ and $0 \leq t_1 < t_2 < \infty$. To see this
note that the first two terms are bounded in $L^2$ by Theorem \ref{velest} whereas the last two terms
were estimated in $L^2$, with a logarithmically growing norm as $\vare \to 0$, in \eqref{fvare}.

Therefore, since $\nabla^{\perp} \eta^{\vare}$ and $\nabla^{\perp}\psi = \mathbb{P}\Phi$ are both square
integrable functions, it follows by the Dominated Convergence Theorem that
\[ \lim_{R\to \infty} A_1 = \int [Eu^{\vare}(t_2,\cdot) - Eu^{\vare}(t_1,\cdot)] \cdot (\nabla^{\perp}\eta^{\vare}) \psi \, dx \]
and
\[ \lim_{R\to \infty} A_2 = \int [Eu^{\vare}(t_2,\cdot) - Eu^{\vare}(t_1,\cdot)] \cdot \eta^{\vare}(\nabla^{\perp}\psi) 
\, dx.\]

Furthermore, it is easy to see that
$\nabla^{\perp} \chi^R$ converges to zero weakly in $L^2$ when $R \to \infty$.
As $[Eu^{\vare}(t_2,\cdot) - Eu^{\vare}(t_1,\cdot)] \cdot \eta^{\vare} \psi$ does not depend on $R$ and
belongs to $L^2$, we infer that $A_3 \to 0$ as $R \to \infty$. 

We have found that the left hand side of identity \eqref{lhs} has a limit as $R\to \infty$. We can rewrite
this limit as follows
\[\lim_{R\to \infty} (A_1 + A_2 + A_3) = \langle \mathbb{P}[\eta^{\vare}Eu^{\vare}](t_2,\cdot) + F^{\vare}(t_2), \Phi\rangle
- \langle \mathbb{P}[\eta^{\vare}Eu^{\vare}](t_1,\cdot) + F^{\vare}(t_1), \Phi\rangle. \]
On the other hand, by identity \eqref{lhs}, and using \eqref{bag1} and \eqref{bag2} we have
\begin{equation} \label{bag3}
\lim_{R\to \infty} |A_1 + A_2 + A_3| \leq \limsup_{R\to \infty} |I_1| + |I_2| \leq C \|\Phi\|_{H^2} |t_2 - t_1|^{1/9},
\end{equation}
which shows that $\mathbb{P}[\eta^{\vare}Eu^{\vare}] + F^{\vare}$ is equicontinuous as a function of time into $H^{-2}$.

We conclude this proof by showing that $\{\mathbb{P}[\eta^{\vare}Eu^{\vare}] + F^{\vare}\}$ is uniformly bounded in 
$L^{\infty}(J;H^{-2}_{\loc}(\real^2))$. We do not need to prove that $\{F^{\vare}\}$ is bounded in this space
because we have already shown that $F^{\vare} \to 0$ as $\vare \to 0$ in $L^{\infty}(J;H^{-2}(\real^2))$.
The only thing left is to prove the boundedness of $\{\mathbb{P}[\eta^{\vare}Eu^{\vare}]\}$.
To this end, let $\chi \in C^{\infty}_c(\real^2)$, fix $p>2$, $1<r<2$ and write
\[|\langle \chi \mathbb{P}[\eta^{\vare}Eu^{\vare}], \Phi \rangle| = |\langle \eta^{\vare}Eu^{\vare}, \mathbb{P}[\chi\Phi]
\rangle| \]
\[\leq  |\langle \eta^{\vare}(Eu^{\vare}- E \mathbf{o}^{\vare}_0 - E \mathbf{i}^{\vare}), \mathbb{P}[\chi\Phi]\rangle| 
+ |\langle \eta^{\vare}E\mathbf{o}^{\vare}_0, \mathbb{P}[\chi\Phi] \rangle| + 
|\langle \eta^{\vare}E\mathbf{i}^{\vare}, \mathbb{P}[\chi\Phi] \rangle| \]
\[\leq \|Eu^{\vare} - E \mathbf{o}^{\vare}_0 - E\mathbf{i}^{\vare}\|_{L^2} \|\mathbb{P}[\chi\Phi]\|_{L^2} + 
\|E\mathbf{o}^{\vare}_0\|_{L^p} \|\mathbb{P}[\chi\Phi]\|_{L^{p/(p-1)}}  + 
\|E\mathbf{i}^{\vare}\|_{L^r} \|\mathbb{P}[\chi\Phi]\|_{L^{r/(r-1)}}\]
\[\leq C(\chi,p,r,J) \|\Phi\|_{H^2}, \]
where in the last inequality we used Theorem \ref{ellpest}, Lemma \ref{idp}, Theorem \ref{velest} and 
relation \eqref{Biotest}. Note that $C$ is independent of $t \in J$. 

It follows from Arzela-Ascoli that, for each $[t_1,t_2] \subset [0,\infty)$ and each $B_R \subset \real^2$, there is a 
subsequence of $\mathbb{P}[\eta^{\vare}Eu^{\vare}]$ which converges strongly in $L^{\infty}([t_1,t_2];H^{-3}(B_R))$. By taking diagonal subsequences we may assume that there is a subsequence which converges strongly in 
$L^{\infty}_{\loc}([0,\infty);H^{-3}_{\loc}(\real^2))$.

This concludes the proof.

\end{proof}

\begin{remark} \label{contH-3}
{\em It follows from the proof of Proposition \ref{tcomp} that any strong limit of
$\{\mathbb{P}[\eta^{\vare}Eu^{\vare}]\}$ in $L^{\infty}_{\loc}([0,\infty);H^{-3}_{\loc}(\real^2))$ 
in fact belongs to $C([0,\infty);H^{-3}_{\loc}(\real^2))$. This is true as 
$\{\mathbb{P}[\eta^{\vare}Eu^{\vare}] + F^{\vare}\}$ is equicontinuous and bounded 
as a function of time into $H^{-2}_{\loc}$ and therefore, by Arzela-Ascoli, its limits are continuous.
Furthermore, $F^{\vare} \to 0$, so that the limits of $\{\mathbb{P}[\eta^{\vare}Eu^{\vare}]\}$ and of
$\{\mathbb{P}[\eta^{\vare}Eu^{\vare}] + F^{\vare}\}$ are the same.}  
\end{remark}

\section{Passing to the limit}

In this section we state and prove our main result. Let us begin with an improvement of the space-time compactness we 
have, which is a consequence of Proposition \ref{tcomp}, obtained by means of interpolation.

\begin{lemma} \label{tcomptoo}
The sequence $\{Eu^{\vare}\}$ is precompact in $L^2_{\loc}((0,\infty)\times \real^2)$.
\end{lemma}

\begin{proof}
By Lemma \ref{sobsilly}, first order derivatives of functions that vanish on 
$\Gamma_{\vare}$ commute with the extension operator, and therefore, for any positive
time, $Eu^{\vare} = \mathbb{P}[Eu^{\vare}]$. We write
\[Eu^{\vare} = \mathbb{P}[(1-\eta^{\vare})Eu^{\vare}] + \mathbb{P}[\eta^{\vare}Eu^{\vare}] \equiv B_1 + B_2. \]
First we note that $B_1 \to 0$ strongly in
$L^2_{\loc}([0,\infty)\times\real^2)$. Indeed, let us fix $0 \leq t_1
< t_2 <\infty$ and set
$J = [t_1,t_2]$. By Theorem \ref{velest}, properties of the cutoff $\eta^{\vare}$ and the fact that the Leray projector
is continuous from $L^2$ to itself we have
\[\|\mathbb{P}[(1-\eta^{\vare})Eu^{\vare}]\|_{L^2(J \times \real^2)} \leq C 
\|(1-\eta^{\vare})Eu^{\vare}\|_{L^2(J \times \real^2)} \] \[\leq C \|1 - \eta^{\vare}\|_{L^4(\real^2)}
\|Eu^{\vare}\|_{L^2(J;L^4(\real^2))} \leq C\vare^{1/2},\]
which proves the desired estimate on $B_1$.

Next we work on $B_2$. We know from Proposition \ref{tcomp} that $B_2$ is precompact in 
$L^{\infty}_{\loc}((0,\infty);H^{-3}_{\loc}(\real^2))$. We will show that, for any $1 < q < 2$,
$B_2$ is bounded in $L^2_{\loc}((0,\infty);W^{1,q}_{\loc}(\real^2))$. The result will follow by 
interpolation. Fix $1<q<2$ and let $q^{\ast} = 2q/(2-q) > 2$.  By Theorem \ref{velest}, $\{Eu^{\vare}\}$
is bounded in $L^4_{\loc}((0,\infty);L^4(\real^2))$. Since $|\eta^{\vare}| \leq 1$ and since $\mathbb{P}$ is
continuous from $L^4(\real^2)$ into itself, it follows that $B_2$ is bounded in 
$L^4_{\loc}((0,\infty);L^4(\real^2))$ which can be continuously imbedded into 
$L^2_{\loc}((0,\infty);L^q_{\loc}(\real^2))$.

What remains is to show that derivatives of $B_2$ are also uniformly
bounded in $L^2_{\loc}((0,\infty);L^q_{\loc}(\real^2))$. Since the
gradient and the Leray projector $\mathbb{P}$ are both Fourier
multipliers, the gradient commutes with $\mathbb{P}$. Therefore,
\[D(\mathbb{P}[\eta^{\vare}E u^{\vare}]) = \mathbb{P}[\eta^{\vare} (D Eu^{\vare})] + 
\mathbb{P}[(D\eta^{\vare}) Eu^{\vare}]\equiv B_{21} + B_{22}. \]

By Theorem \ref{velest}, $DEu^{\vare}$  is bounded in $L^2_{\loc}((0,\infty);L^2(\real^2))$.
Since $|\eta^{\vare}| \leq 1$ and $\mathbb{P}$ is continuous from $L^2(\real^2)$ to itself, we immediately
obtain the desired estimate for $B_{21}$. As for the term $B_{22}$, we use Theorem \ref{velest} once again
to obtain that $Eu^{\vare}$ is bounded in  $L^{2q^{\ast}/(q^{\ast}-2)}_{\loc}((0,\infty);L^{q^{\ast}}(\real^2))$ 
and we recall that $D\eta^{\vare}$ is uniformly bounded in $L^2(\real^2)$. With this, we have that
$(D\eta^{\vare}) Eu^{\vare}$  is bounded in $L^{2q^{\ast}/(q^{\ast}-2)}_{\loc}((0,\infty);L^q(\real^2))$,
continuously imbedded into $L^2_{\loc}((0,\infty);L^q(\real^2))$. This concludes the proof.

\end{proof}

We will prove that limits of the sequence $\{Eu^{\vare}\}$ are solutions of the Navier-Stokes equations
in a suitable weak sense. To be precise, we formulate the notion of weak solution we will use.

\begin{definition} \label{wform}
Let $u \in L^2_{\loc}((0,\infty)\times \real^2)\cap C([0,\infty);\mathcal{D}^{\prime}(\real^2))$. 
We say that $u$ is a weak solution of the incompressible
Navier-Stokes equations with initial velocity $u_0$ if, for any divergence-free test vector field 
$\Psi \in C^{\infty}_c((0,\infty)\times \real^2)$, we have
\[\int_0^{\infty} \int_{\real^2} \bigl(u \cdot \Psi_t +  [(u \cdot \nabla) \Psi] \cdot u  + 
\nu u \cdot \Delta \Psi\bigr) \, dxdt = 0.\]
Furthermore, for every $t \geq 0$, $\mbox{div }u(t,\cdot) = 0$ in the sense of distributions and
$u(t,\cdot) \rightharpoonup u_0$ in the sense of distributions as $t \to 0^+$. 
\end{definition}

Recall that $\mathcal{K}$ denotes the kernel of the Biot-Savart law, as introduced in \eqref{biotker}. We are finally
ready to state and prove the main result of this work.

\begin{theorem} \label{mainthm}
Any strong limit $u$ of $\{Eu^{\vare}\}$ in $L^2_{\loc}((0,\infty)\times \real^2)$ is 
a weak solution of the incompressible Navier-Stokes equations in $\real^2$ with initial 
velocity given by $u_0 = \mathcal{K} \ast\omega_0 + \gamma H$. 
\end{theorem}

\begin{remark} \label{exist} {\em As $\{Eu^{\vare}\}$ is precompact, by virtue of Lemma \ref{tcomptoo}, 
there exists at least one such strong limit.}
\end{remark}

\begin{proof}

For each $\vare$ sufficiently small, choose $0 < \delta < 1$ such that $\{|x| > 2 \delta\} \subseteq \Pi_{\vare}$.
Clearly, if $\{|x| > 2 \delta\} \subseteq \Pi_{\vare_0}$  then $\{|x| > 2 \delta\} \subseteq \Pi_{\vare}$,
for all $\vare \leq \vare_0$. Also consider $R > 2 > 2\delta$. We use the cutoff $\varphi$ introduced in
Section 4 to define: 
\[\varphi^{\delta} = \varphi^{\delta}(x) \equiv \varphi(|x|/\delta) \mbox{ and } \chi^R = \chi^R(x) \equiv 
1 - \varphi(|x|/R). \]

As in the proof of Proposition \ref{tcomp} we let $\Phi$ be a smooth, compactly supported vector 
field in $\real^2$, which, in addition, we assume to be divergence-free. 
We define $\psi = \psi(x) = [(\Delta)^{-1} \mbox{ curl } \Phi](x) -  
[(\Delta)^{-1} \mbox{ curl } \Phi](0)$. Recall that $\nabla^{\perp}\psi =  \Phi$
and that $\psi$ satisfies \eqref{veryclever}.
We also consider $\theta = \theta(t) \in C^{\infty}_c((0,\infty))$. 

We use the test function $\varphi^{\delta} \theta \psi \chi^R$, which belongs to 
$C^{\infty}_c((0,\infty) \times \Pi_{\vare})$ in the weak form of the vorticity equation.
We can rewrite the integrals on $\Pi_{\vare}$ as full plane integrals using the extension operator
to obtain the following integral identity
\[\int_0^{\infty} \int_{\real^2} E\omega^{\vare} \theta_t \varphi^{\delta} \psi \chi^R \,dxdt 
+ \int_0^{\infty} \int_{\real^2} Eu^{\vare} E\omega^{\vare} \cdot \theta \nabla(\varphi^{\delta} \psi \chi^R) \,dxdt\]
\begin{equation} \label{grab} + \nu \int_0^{\infty} \int_{\real^2} E\omega^{\vare} \theta \Delta(\varphi^{\delta} \psi \chi^R) \,dxdt = 0. 
\end{equation} 

Our first step is to pass to the limit $\vare \to 0$ in this identity, while keeping $\delta$ and $R$ fixed.
Let $u$ be a strong limit in $L^2_{\loc}((0,\infty)\times \real^2)$ of a subsequence $Eu^{\vare_k}$ 
of $Eu^{\vare}$. We observe that $E\omega^{\vare_k} \to \mbox{ curl }u \equiv \omega$ 
strongly in $L^2_{\loc}((0,\infty);H^{-1}_{\loc}(\real^2))$. Similarly, 
we may also deduce that $u$ is divergence-free in the sense of distributions. The passage to the limit is 
immediate in the linear terms of \eqref{grab}. For the nonlinear term we recall that, by Theorem \ref{velest}, 
$\{E\omega^{\vare_k}\}$ is uniformly bounded in $L^2_{\loc}((0,\infty);L^2(\real^2))$. Hence a subsequence of
$\{E\omega^{\vare_k}\}$ converges weakly in $L^2_{\loc}((0,\infty);L^2(\real^2))$. Using the convergence 
$E\omega^{\vare_k} \to \omega$ strong in $L^2_{\loc}((0,\infty);H^{-1}_{\loc}(\real^2))$ and uniqueness of
weak limits we conclude that $E\omega^{\vare_k} \rightharpoonup \omega$ weakly in 
$L^2_{\loc}((0,\infty);L^2(\real^2))$,
without passing to further subsequences. Now $Eu^{\vare_k} E\omega^{\vare_k}$ is a weak-strong pair, so that
we can pass to the limit in the nonlinear term as well. We arrive at the identity 
\[J_1 + J_2 + J_3 \equiv \int_0^{\infty} \int_{\real^2} \omega \,\theta_t \varphi^{\delta} \psi \chi^R \,dxdt 
+ \int_0^{\infty} \int_{\real^2} u\,\omega \,\cdot \theta \nabla(\varphi^{\delta} \psi \chi^R) \,dxdt\]
\begin{equation} \label{varelim}
+ \, \nu \int_0^{\infty} \int_{\real^2} \omega \,\theta \Delta(\varphi^{\delta} \psi \chi^R) \,dxdt = 0. 
\end{equation} 

Now we pass to the limit both $\delta \to 0$ and $R \to \infty$ in each separate term in \eqref{varelim}.
We begin with $J_1$.

First we observe that
 \begin{equation} \label{obs}
u - \alpha \varphi \left( \frac{\beta|x|}{\lambda_0}\right) H  \equiv u - F \in L^{\infty}_{\loc}((0,\infty);L^2(\real^2)),
\end{equation}
where $\alpha$ was introduced in \eqref{u0}, $\beta$ in Lemma \ref{confmap} and $\lambda_0$ in Lemma \ref{idp}
and $H$ is the harmonic vector field introduced in \eqref{helena}.
Identity \eqref{obs} follows from the convergence $Eu^{\vare_k} \to u$, from the fact that $Eu^{\vare} - E\mathbf{o}_0^{\vare}$ 
is bounded in $L^{\infty}_{\loc}((0,\infty);L^2(\real^2))$, see Theorem \ref{velest} and from the fact that 
$E\mathbf{o}_0^{\vare} \to \alpha \varphi \left( \frac{\beta|x|}{\lambda_0}\right) H$ strongly in $L^2(\real^2)$
by Lemma \ref{idp}. 

Next, write $\omega = \nabla^{\perp} \cdot u$ and integrate by parts to obtain
\[ J_1 = - \int_0^{\infty} \int_{\real^2} u \cdot \theta_t \nabla^{\perp} (\varphi^{\delta} \psi \chi^R) \,dxdt 
= - \int_0^{\infty} \int_{\real^2} (u-F)  \cdot \theta_t \nabla^{\perp} (\varphi^{\delta} \psi \chi^R) \,dxdt, \]
where we have used the fact that $F$ does not depend on time, so that the additional integral vanishes.

We write
\[ J_1 = - \int_0^{\infty} \int_{\real^2} (u-F)  \cdot \theta_t (\nabla^{\perp}\varphi^{\delta}) \psi \chi^R \,dxdt
- \int_0^{\infty} \int_{\real^2} (u-F)  \cdot \theta_t \varphi^{\delta}  \Phi \chi^R \,dxdt \]
\[ - \int_0^{\infty} \int_{\real^2} (u-F)  \cdot \theta_t \varphi^{\delta} \psi (\nabla^{\perp} \chi^R) \,dxdt 
\equiv - J_{11} - J_{12} - J_{13}. \]

It is easy to see that $\nabla \varphi^{\delta}$ converges to zero weakly in $L^2(\real^2)$ 
when $\delta\to 0$ and therefore, 
\begin{equation} \label{coitado11}
\lim_{R\to \infty} \lim_{\delta \to 0} J_{11} = 0.
\end{equation}

On the other hand,  
$\nabla \chi^R$ also converges to zero weakly in $L^2(\real^2)$ when $R\to \infty$. Furthermore, 
$\varphi^{\delta} \to 1$ poinwise when $\delta \to 0$, so by dominated convergence, we find that
\begin{equation} \label{coitado12}
\lim_{R\to \infty} \lim_{\delta \to 0} J_{13} = 0.
\end{equation}

We also note that $\varphi^{\delta} \chi^R$ converges poinwise to $1$ as $\delta \to 0$ and $R \to \infty$
(no matter which order), so that, by dominated convergence, we deduce that
\begin{equation} \label{coitado13} \lim_{R \to \infty} \lim_{\delta \to 0} J_{12} 
= \int_0^{\infty} \int_{\real^2} (u-F)  \cdot \theta_t  \Phi\,dxdt = \int_0^{\infty} \int_{\real^2} u  \cdot \theta_t  \Phi\,dxdt.
\end{equation}

Putting together \eqref{coitado11}, \eqref{coitado12} and \eqref{coitado13} we obtain
\begin{equation} \label{coitado1} 
\lim_{R\to\infty}\lim_{\delta\to 0} J_1 = - \int_0^{\infty} \int_{\real^2} u \cdot \theta_t  \Phi\,dxdt .
\end{equation}

Next we treat the nonlinear term $J_2$. First note that the uniform estimates on $Eu^{\vare}$ contained
in Theorem \ref{velest} imply that $u \in L^4_{\loc}((0,\infty);L^4(\real^2))$. The argument is the same
we used to prove $\omega \in L^2_{\loc}((0,\infty);L^2(\real^2))$

We write
\[J_2 = \int_0^{\infty} \int_{\real^2} u\,\omega \,\cdot \theta (\nabla\varphi^{\delta}) \psi \chi^R \,dxdt
-\int_0^{\infty} \int_{\real^2} u\,\omega \,\cdot \theta \varphi^{\delta} \Phi^{\perp} \chi^R \,dxdt \]
\[+ \int_0^{\infty} \int_{\real^2} u\,\omega \,\cdot \theta \varphi^{\delta} \psi (\nabla\chi^R) \,dxdt \equiv
J_{21} - J_{22} + J_{23}. \]

We have that
\[|J_{21}| \leq \int_0^{\infty} |\theta| \|u\|_{L^4} \|\omega\|_{L^2} \|\psi\nabla\varphi^{\delta}\|_{L^4} \, dt 
= \mathcal{O}(\sqrt{\delta}) \]
and
\[|J_{23}| \leq \|\psi\|_{L^{\infty}}\int_0^{\infty} |\theta| \|u\|_{L^4} \|\omega\|_{L^2} \|\nabla\chi^R\|_{L^4} \, dt 
= \mathcal{O}(R^{-1/2}).\]

We conclude that
\begin{equation} \label{coitado21}
\lim_{R \to \infty} \lim_{\delta \to 0} (J_{21} + J_{23}) = 0.
\end{equation}

In addition, by dominated convergence  we have that  
\begin{equation} \label{coitado22}
\lim_{R\to \infty} \lim_{\delta \to 0} J_{22} =  \int_0^{\infty} \int_{\real^2} u\,\omega \,\cdot \theta \Phi^{\perp} \,dxdt = - \int_0^{\infty} \int_{\real^2} (u\cdot\nabla)u \theta \cdot \Phi \,dxdt,
\end{equation}
where this last equality follows from the identity $u\cdot\nabla u - (u\omega)^{\perp} = \nabla(|u|^2/2)$,
together with the fact that $\Phi$ is divergence free.

Therefore, using \eqref{coitado21} and \eqref{coitado22} and integrating by parts we find
\begin{equation} \label{coitado2}
\lim_{R\to\infty}\lim_{\delta\to 0} J_2 = - \int_0^{\infty} \int_{\real^2} [(u\cdot\nabla) \theta\Phi] \cdot u \,dxdt.
\end{equation}

Lastly we treat $J_3$. Once again, we write
\[J_3 = \nu \int_0^{\infty} \int_{\real^2} \omega \,\theta (\Delta\varphi^{\delta}) \psi \chi^R \,dxdt
 + \nu \int_0^{\infty} \int_{\real^2} \omega \,\theta \varphi^{\delta} (\Delta\psi) \chi^R \,dxdt \]
\[+ \nu \int_0^{\infty} \int_{\real^2} \omega \,\theta \varphi^{\delta} \psi (\Delta \chi^R) \,dxdt
+ 2 \nu \int_0^{\infty} \int_{\real^2} \omega \,\theta \nabla\psi \cdot ((\nabla \varphi^{\delta})\chi^R
+\varphi^{\delta}(\nabla\chi^R)) \,dxdt\] \[\equiv J_{31} + J_{32} + J_{33} + J_{34}.\]

Using, similarly to what we have already done, that: $\omega \in L^2_{\loc}((0,\infty);L^2(\real^2))$, $\psi \Delta \varphi^{\delta} \rightharpoonup 0$ weakly in $L^2(\real^2)$ as $\delta \to 0$, $\Delta \chi^R \to 0$ strongly in $L^2(\real^2)$ as $R \to \infty$, $\nabla \varphi^{\delta} \rightharpoonup 0$  weakly in $L^2(\real^2)$ as $\delta \to 0$, $\nabla \chi^R \rightharpoonup 0$  weakly in $L^2(\real^2)$ as $R \to \infty$ and dominated convergence, we deduce that    
\begin{equation} \label{coitado31}
\lim_{R\to \infty}\lim_{\delta\to 0} (J_{31}+J_{33}+J_{34}) = 0.
\end{equation}
Therefore we obtain, integrating by parts, 
\begin{equation} \label{coitado3}
\lim_{R\to\infty}\lim_{\delta \to 0} J_3 = 
\nu \int_0^{\infty} \int_{\real^2} \omega \,\theta (\Delta\psi) \,dxdt = - 
\nu \int_0^{\infty} \int_{\real^2} u \cdot \,\theta \Delta \Phi \,dxdt. 
\end{equation}

Recall that $J_1 + J_2 + J_3 = 0$, so that, adding \eqref{coitado1} with \eqref{coitado2} and \eqref{coitado3} we find
\begin{equation} \label{limitpde}
\int_0^{\infty}\int_{\real^2} \bigl(u \cdot \theta_t \Phi + [(u\cdot\nabla) \theta \Phi]\cdot u + \nu u\cdot \theta \Delta \Phi\bigr)  \, dxdt = 0.
\end{equation}
We observe that linear combinations of products of smooth, compactly
supported functions of the form $\theta \,\Phi$ are dense in
$C^{\infty}_c((0,\infty)\times \real^2)$. With this observation and
\eqref{limitpde} we find that $u$ satisfies the integral identity in
Definition \ref{wform}.  We have already noted that $u$ is
divergence-free in the sense of distributions. All that remains is to show that
$u \in C([0,\infty);\mathcal{D}^{\prime}(\real^2))$ and that
$u(t,\cdot) \rightharpoonup u_0$ in $\mathcal{D}^{\prime}$ as $t \to 0$.

First note that in the proof of Lemma \ref{tcomptoo} we showed that 
$Eu^{\vare} - \mathbb{P}[\eta^{\vare}Eu^{\vare}] \to 0$ strongly in
$L^2_{\loc}((0,\infty) \times \real^2)$. 
Therefore, $\mathbb{P}[\eta^{\vare_k}Eu^{\vare_k}] \rightarrow u$ in
$L^2_{\loc}((0,\infty) \times \real^2)$. 
By Remark \ref{contH-3}, we have 
a subsequence of $\{\mathbb{P}[\eta^{\vare_k}Eu^{\vare_k}]\}$ which
converges strongly in  
$L^{\infty}_{\loc}([0,\infty);H^{-3}_{\loc}(\real^2))$ to a limit $v
\in C([0,\infty);H^{-3}_{\loc}(\real^2))$. 
It follows by uniqueness of limits (in
$L^2_{\loc}((0,\infty);H^{-3}_{\loc})$, for example)  
that $u(t,\cdot)=v(t,\cdot)$ for almost all $t \in (0,\infty)$, which 
implies 
that $u$ can be identified with $v$. This in turn implies that $u \in 
C([0,\infty);\mathcal{D}^{\prime}(\real^2))$. Furthermore, as
$\mathbb{P}[\eta^{\vare_k}Eu^{\vare_k}]$ converges to $u$ uniformly in
time with values in  $H^{-3}_{\loc}$, one has that
$\mathbb{P}[\eta^{\vare_k}Eu^{\vare_k}_0]$ converges to $u_0$ in
$H^{-3}_{\loc}$. On the other hand, Lemma \ref{limidp} says that 
$\mathbb{P}[\eta^{\vare_k}Eu^{\vare_k}_0]$ converges to
$\mathcal{K}\ast\omega_0+\gamma H$ in the sense of distributions. By
uniqueness of the limit in $\mathcal{D}'$, we conclude that
$u_0=\mathcal{K}\ast\omega_0+\gamma H$. 
This  concludes this proof.
\end{proof}

\begin{remark}\label{initdatacont}
{\em At the end of the proof above we showed that the initial data for the limit problem is attained in $C([0,\infty);\mathcal{D}^{\prime}(\real^2))$. We can actually prove a stronger statement, namely that 
there exists a positive constant $C>0$ such that 
\[|\langle u(t)-u_0, \phi \rangle | \leq C\|\phi\|_{H^2}t^{1/9},\]
for all $\phi \in \mathcal{D}$. Indeed, to see this fix $\phi \in C^{\infty}_c(\real^2)$ and let $\psi = \psi(x) = [\Delta^{-1}\mbox{ curl }\phi](x) - [\Delta^{-1}\mbox{ curl }\phi](0)$. Consider the sequence of approximations $\{u^{\vare_k}\}$ constructed in the proof above. For each $\vare_k$ recall the 
auxiliary functional $F^{\vare_k}$ used in the proof of Proposition \ref{tcomp}, given by
\[ \langle F^{\vare_k}(t) , \phi \rangle  =  \int (Eu^{\vare_k}(t,x) - E\mathbf{o}_0^{\vare_k}(x)) \cdot 
(\nabla^{\perp}\eta^{\vare_k})(x) \psi(x) \,dx. \]  
Write
\begin{equation} \label{lall}
\langle u(t) - u_0 , \phi \rangle   \equiv  L_1 + L_2 + L_3 + L_4,
\end{equation}
where
\[ L_1 = \langle u(t) - \mathbb{P}[\eta^{\vare_k} E u^{\vare_k}](t) - F^{\vare_{k}}(t) , \phi \rangle;\]
\[ L_2 = \langle \mathbb{P}[\eta^{\vare_k} E u^{\vare_k}](t) + F^{\vare_{k}}(t)  - \mathbb{P}[\eta^{\vare_k} E u^{\vare_k}_0] - F^{\vare_{k}}(0), \phi \rangle; \]
\[ L_3 = \langle F^{\vare_{k}}(0), \phi \rangle; \]
\[ L_4 = \langle \mathbb{P}[\eta^{\vare_k} E u^{\vare_k}_0] - u_0, \phi \rangle.\]

By Remark \ref{contH-3}, we have that $\lim_{\vare_k \to 0} L_1 = 0$ uniformly in time.
By the argument in the proof of Proposition \ref{tcomp}, see estimate \eqref{bag3}, we find that
$|L_2| \leq C \|\phi\|_{H^2} t^{1/9}$. By estimate \eqref{fvare} we know that $\lim_{\vare_k \to 0} L_3 = 0$.
Finally, by Lemma \ref{limidp} we obtain that  $\lim_{\vare_k \to 0} L_4 = 0$. Therefore, using \eqref{lall}
we deduce
\[ |\langle u(t) - u_0 , \phi \rangle| \leq C \|\phi\|_{H^2} t^{1/9}, \]
as desired.}
\end{remark}

\begin{remark}\label{regul}
{\em The solution obtained at the limit is much more regular then what is
required for a weak solution as stated in Definition \ref{wform}. For
example, the a priori estimates given in Theorem \ref{velest} imply
immediately that $u\in L^{\frac{2p}{p-2}}_{\loc}((0,\infty);L^p)$ for
all $p\in(2,\infty)$ and $\nabla u\in   L^{2}_{\loc}((0,\infty);L^2)$.} 
\end{remark}

Our result above provides strong compactness of viscous flows around a small obstacle but does not
address actual convergence. The passage from compactness to convergence is clearly reduced
to the issue of uniqueness for the limit problem.  The issue of uniqueness of solutions for the 2D 
incompressible Navier-Stokes equations with initial data of the form $u_0$ is classical and delicate.
Let us briefly review the related  literature.  The first relevant results are due to G. Benfatto, 
R. Esposito and M. Pulvirenti, see \cite{BEP}, who showed uniqueness for initial flows of the form
$\sum \gamma_i H(x-x_i)$ if $\sum|\gamma_i|$ is sufficiently small,
and to G.-H. Cottet, see \cite{Cottet}, who showed uniqueness for
small initial vorticities which are general bounded measures.  Later, 
Y. Giga, T. Miyakawa and H. Osada generalized this uniqueness result for initial flows of the form
$\mathcal{K}\ast\omega$ with $\omega$ a Radon measure with sufficiently small atomic part, see \cite{GMO,kato94}.  
As we observed in the Introduction, this smallness condition is closely related, in a technical sense, 
to the smallness condition on $\gamma$ which we also had to impose, see \eqref{conditions}. 
Recently the uniqueness assumption on the 
atomic part of $\omega$ was removed, first by T. Gallay and C. E. Wayne for initial flow of the form 
$\gamma H$ in \cite{GW05}, see also \cite{GGL}, and then for general 
$\mathcal{K}\ast\omega_0$ initial flows with $\omega_0$ an arbitrary Radon measure by I. Gallagher and 
T. Gallay in \cite{GG05}. These results are a byproduct of the remarkable large-time asymptotics 
results obtained by Gallay and Wayne in \cite{GW02}.

From the point of view of the present work, the natural question is whether Gallagher and Gallay's 
result in \cite{GG05} implies uniqueness of the weak solutions obtained as small obstacle limits
in Theorem \ref{mainthm}. Gallagher and Gallay's proved that there is at most one solution $v = v(x,t)$ 
to the 2D Navier-Stokes equations with initial vorticity a Radon measure satisfying
\begin{enumerate}
\item the vorticity $w = \mbox{curl }v \in C((0,T);L^1(\real^2) \cap L^{\infty}(\real^2))$.
\item $\|w(\cdot,t\|_{L^1} \leq K$ for all $t>0$.
\item $w(\cdot,t) \rightharpoonup \omega_0 + \gamma \delta_0$ as $t \to 0$
\item the vorticity $w$ is a mild solution of the vorticity formulation of Navier-Stokes equations on any time interval $[t_0,t_1]$ compactly contained in $(0,\infty)$, i.e.,
\[w(t) = e^{\nu (t-t_0)\Delta}w(t_0) - \int_{t_0}^t \nabla \cdot e^{\nu (t-s) \Delta} [v(s)w(s)]\, ds, \]
for any $0 < t_0 < t < t_1$.
\end{enumerate}

In our result, we obtain a weak solution $u$ in the sense of Definition 1 with initial vorticity 
$\omega_0 + \gamma \delta$. The additional regularity stated in Remark
\ref{regul} should be more than enough to prove that the solution $u$
is a strong solution for $t>0$, thus placing our result in the setting
of Gallagher and Gallay. For example, the condition 
$u \in L^{q_0}_{\loc}((0,\infty);L^p(\real^2))$ places our solution
on the Serrin criterion curve, 
see \cite{serrin62}, which is a standard condition for uniqueness of weak solutions, see also 
\cite{lemarie02}. 

In short, there is no uniqueness result stated in the literature which includes precisely
our solution, but clearly such a result is expected and should be easy
to prove. We will not address this issue any further, as it 
escapes the main purpose of this paper.           

\section{Conclusions}
   The purpose of this section is to interpret what we have done  
in a broader context and to point out some directions for improvement and further work. Our basic 
problem was to find conditions under which the presence of a single small obstacle could be ignored 
in the modelling of large scale flow. The precise formulation we used, working in the unbounded exterior domain and 
fixing the large scale flow by choosing an initial vorticity $\omega_0$ and a 
circulation $\gamma$, was convenient from the mathematical point of view, but it was far from physically natural. 

     There are many ways in which to formulate mathematically the problem of placing a small 
obstacle in a given incompressible viscous flow. To be more precise, let us suppose that
we are given a smooth background flow $u_b = u_b(t,x)$ and say we wish to insert a small circular 
obstacle (centered at the origin) within this flow. That would mean adding a correction 
$u_c^{\vare} = u_c^{\vare}(t,x)$ so that the new flow $u^{\vare} \equiv u_b + u_c^{\vare}$ satisfies the 
no-slip boundary condition at $|x|=\vare$. It makes sense to assume, for simplicity, that there are no more 
boundaries in the problem and that $u_b$ vanishes at infinity. We also expect $u_c^{\vare}$ to be sharply 
localized, and hence it should vanish at infinity as well. In this context, there is no reason for the background 
vorticity $\omega_b = \mbox{ curl }u_b$ to have vanishing total integral. Let us denote by $m_b$ the total 
integral of $\omega_b$. The vanishing of $m_b$ is equivalent to requiring the background 
flow to have globally finite kinetic energy, a physically reasonable assumption. However, 
infinite energy flows associated with smooth but mass-unbalanced vorticity have often been 
considered (see Section 3.1.3 in \cite{MB02} for a discussion), so we do not assume $m_b = 0$, but we consider 
this case to be of special physical relevance.  

Let $\omega^{\vare} = \mbox{ curl } u^{\vare}$. As we have seen, there exists $\alpha = \alpha(\vare)$ such that 
$u^{\vare} = K^{\vare}[\omega^{\vare}] + \alpha(\vare) H$, where $H$ was defined in \eqref{helena}. 
If $u^{\vare}$ should converge to $u_b$ then $\int \omega^{\vare} - \alpha(\vare) \to 0$. If this
were not the case, then the kinetic energy of $u^{\vare}$ would blow up near the origin, 
a complicated consequence of the method of images for the disk, see \cite{iln-shrink,saffman}.  
Moreover, since both $|u_b  - m_b H| = \mathcal{O}(1/|x|^2)$ and $|u^{\vare} - \alpha(\vare) H| 
= \mathcal{O}(1/|x|^2)$ at infinity (see \cite{MB02}), we must have $\alpha(\vare)  \to m_b$ as $\vare \to 0$
for the kinetic energy of the difference $u_b - u^{\vare}$ to be bounded at infinity. 

In other words, at any fixed time, if a flow $u^{\vare}$ around the obstacle approximates the background
flow $u_b$ without infinite energy discrepancies either at infinity or locally, then it must have vorticity
$\omega^{\vare}$ with total integral close to $m_b$ and its harmonic part should also be close to $m_b$. 
It is therefore very natural to model this problem with a solution of the Navier-Stokes equations for which:
\begin{enumerate}
\item the initial vorticity $\omega_0^{\vare}$ has integral approaching $m_b$,
\item $\omega_0^{\vare}$ approximates $\omega_b(0,x)$ in $L^1$,
\item the initial circulation $\gamma(\vare) = \alpha(\vare) - \int \omega_0^{\vare} \to 0$ as $\vare \to 0$. 
\end{enumerate}
The main point of this discussion is to claim that the results obtained in our 
article can be easily adapted to include physically natural approximations such as those described above. 
Therefore, although we chose a very specific way of formulating the small obstacle limit, 
reasonable alternatives should lead to the same result. 

From the discussion above and the physical interpretation of the case $\gamma(\vare) \neq  0$, given in Section 5,
together with the natural scaling of this problem, we can see that the cases 
$\gamma (\vare)= 0$ or $\gamma(\vare) \to 0$ are by far the most interesting situations. However, there would have been no substantial simplification of the argument by restricting our problem to the case $\gamma(\vare) = o(1)$. Moreover, the situation in which $\gamma(\vare) = \mathcal{O}(1)$ is mathematically very interesting. 
Indeed, the smallness condition on the initial circulation only appears when we assume $\gamma(\vare) = \mathcal{O}(1)$. Also, there is a discrepancy between the results obtained for the inviscid and viscous cases when $\gamma(\vare) = \mathcal{O}(1)$, which suggests that the limits $\vare \to 0$ and $\nu \to 0$ do not commute. 
         
We have assumed throughout that $\omega_0$ was smooth and compactly supported. How much regularity on
$\omega_0$ did we really use? The answer is none. We actually needed $u_0$ bounded in $L^2_{\loc}$ and
$L^{2,\infty}$ and nothing else. We contrast this with the inviscid argument, where we needed 
$\omega_0$ in $L^p$, $p>2$.
    
Let us turn to some problems which arise naturally from our work.
One particularly interesting question is the issue of considering both the viscosity and the
obstacle small. This should be a difficult problem, because the wake due to an obstacle 
becomes more pronounced and turbulent as viscosity vanishes. It is well-known that, for full plane flow, 
one can take the vanishing viscosity limit, obtaining solutions of the incompressible 
2D Euler equations, see \cite{chemin96,kelliher04,majda93}. In the presence of a material boundary, the vanishing viscosity limit
is a classical open problem, even if the flow is very smooth. The difficulty is due to the boundary layer.
The problem of taking the vanishing viscosity limit outside a very small obstacle interpolates nicely between 
the full plane result and the open problem of taking the vanishing viscosity limit in the presence of a fixed material boundary.
In fact, this question is one of the main motivations of the present work and it is still under 
consideration by the authors. Taking into account the result obtained in this paper it is clear that one should first pursue the small viscosity problem in the case $\gamma(\vare) = 0$ or $\gamma(\vare) = o(1)$ as $\vare \to 0$, since the smallness condition in our convergence result gets more
restrictive as viscosity vanishes. With this future work in mind we have included the specific dependence of our estimates on viscosity for as long as it was practical. 

A second problem is to extend our analysis to velocity fields which are constant at infinity, in order 
to include the classical case of a material body moving in a fluid with roughly constant speed. 

Yet another problem that arises from our work is to remove the smallness condition on the initial circulation. The parallel
between our convergence problem and uniqueness for the limit flow suggests a strategy.
Is it possible to adapt the entropy-entropy flux techniques used by Gallay and Gallagher for the uniqueness
problem to the small obstacle asymptotics? 

A fourth problem is to obtain an asymptotic description of the correction term in the small obstacle limit, i.e. a description of the ``wake" associated with the small obstacle. Finally, one can consider a whole host of related problems, described loosely as the study of limit flows in singularly perturbed domains. For instance, one can  study limit flows in a bounded domain with  one or more small obstacles, or in a domain composed of a small neck joining two fat domains, or in a domain having a long thin tail, etc.
 
\vspace{0.5cm}

\small{ {\it Acknowledgments:} The authors wish to thank Russel Caflisch for a helpful 
discussion. D. I. would like to thank the hospitality of IMECC-UNICAMP, while    
M. C. L. F and H. J. N. L. wish to thank the hospitality of
the Department of Mathematics of the Univ. de Rennes 1. M. C. L. F.'s research is supported 
in part by CNPq grant \# 302.102/2004-3. H. J. N. L.'s research is supported in part by CNPq grant 
\# 302.214/2004-6. This work was supported in part by FAPESP Grant \# 02/02370-5,
by FAEPEX, by the Brazil-France Cooperation and by the Differential Equations PRONEX 
at UNICAMP.}

 \noindent
{\sc Drago\c{s} Iftimie \\
Institut Camille Jordan\\
Universit{\'e} Claude Bernard Lyon 1\\
B{\^a}timent Braconnier (ex-101)\\
21 Avenue Claude Bernard\\
69622 VILLEURBANNE Cedex, France
\\}
{\it E-mail address:} dragos.iftimie@univ-lyon1.fr

\vspace{.1in}

\noindent
{\sc
Milton C. Lopes Filho\\
Departamento de Matematica, IMECC-UNICAMP.\\
Caixa Postal 6065, Campinas, SP 13081-970, Brasil
\\}
{\it E-mail address:} mlopes@ime.unicamp.br

\vspace{.1in}

\noindent
{\sc Helena J. Nussenzveig Lopes\\
Departamento de Matematica, IMECC-UNICAMP.\\
Caixa Postal 6065, Campinas, SP 13081-970, Brasil
\\}
{\it E-mail address:} hlopes@ime.unicamp.br


\begin{thebibliography}{99}
\bibitem{batchelor67} G. K. Batchelor. \textit{An introduction to fluid dynamics.} Second paperback edition, Cambridge Mathematical Library, Cambridge University Press, Cambridge, 1999.

\bibitem{BEP} G. Benfatto; R. Esposito; M. Pulvirenti. \textit{Planar Navier-Stokes flow for singular initial data.} Nonlinear Anal. {\bf 9} (1985), no. 6, 533--545.

\bibitem{chemin96} J.-Y. Chemin. \textit{A remark on the inviscid limit for two-dimensional incompressible fluids.} 
Comm. Part. Diff. Eqs. {\bf 21} (1996) 1771--1779.

\bibitem{Cottet} G.-H. Cottet. \textit{Équations de Navier-Stokes dans le plan avec tourbillon initial mesure.} 
C. R. Acad. Sci. Paris Sér. I Math. {\bf 303} (1986), no. 4, 105--108.

\bibitem{DS1} Wakako Dan; Yoshihiro Shibata. \textit{On the $L\sb q$--$L\sb r$ estimates of the Stokes semigroup in a two-dimensional exterior domain.} J. Math. Soc. Japan {\bf 51} (1999), no. 1, 181--207. 

\bibitem{DS2} Wakako Dan; Yoshihiro Shibata.  \textit{Remark on the $L\sb q$-$L\sb \infty$ estimate of the Stokes semigroup in a $2$-dimensional exterior domain.} Pacific J. Math. {\bf 189} (1999), no. 2, 223--239.

\bibitem{FK} T. Kato; H. Fujita. \textit{On the nonstationary Navier-Stokes
system.} Rend. Sem. Mat. Univ. Padova {\bf 32} (1962), 243--260. 

\bibitem{kelliher04} J. Kelliher. \textit{The inviscid limit for two-dimensional incompressible fluids with 
unbounded vorticity.} Math. Res. Lett. {\bf 11} (2004) 519--528.

\bibitem{GG05} I. Gallagher; T. Gallay. {\it Uniqueness for the two-dimensional Navier-Stokes equation with
a measure as initial vorticity.} Math. Ann. {\bf 332} (2005), 287--327.

\bibitem{GGL} I. Gallagher; T, Gallay; P.-L. Lions. {\it On the uniqueness of the solution of the
two-dimensional Navier-Stokes equation with a Dirac mass as initial vorticity.} ArXiV preprint 
math.AP/0410344, 2004. 

\bibitem{GW02}  T. Gallay; C. E. Wayne. {\it Invariant manifolds and the long-time asymptotics of the Navier-Stokes and vorticity equations on $\real^2$.} Arch. Ration. Mech. Anal. {\bf 163} (2002), no. 3, 209--258.

\bibitem{GW05}  T. Gallay; C. E. Wayne. {\it Global stability of vortex solutions of the two-dimensional 
Navier-Stokes equations.} Comm. Math. Phys. {\bf 255} (2005), 97--129.

\bibitem{giga81} Y. Giga. {\it Analyticity of the semigroup generated by the Stokes operator in $L_r$ spaces.}
Math. Z. {\bf 178} (1981), 297--329.

\bibitem{GMO} Y. Giga; T. Miyakawa; H. Osada. {\it  Two-dimensional Navier-Stokes flow with measures as initial vorticity.} Arch. Rational Mech. Anal. {\bf 104} (1988), no. 3, 223--250.

\bibitem{iln-shrink} D. Iftimie; M. C. Lopes Filho; H. J. Nussenzveig Lopes. {\it Two dimensional incompressible ideal flow around a small obstacle.} Comm. Part. Diff. Eqs. {\bf 28} (2003) 349--379.


\bibitem{kato94} T. Kato.  {\it The Navier-Stokes equation for an incompressible fluid in $\real^2$ with a measure as the initial vorticity.} Differential Integral Equations {\bf 7} (1994), no. 3-4, 949--966.

\bibitem{kozono-yamazaki} H. Kozono; M. Yamazaki. \textit{Local and global unique solvability of the Navier-Stokes exterior problem with Cauchy data in the space $L\sp {n,\infty}$.} Houston J. Math. 21 (1995), no. 4, 755--799. 

\bibitem{lady69} O. Ladyzhenskaya. \textit{The Mathematical Theory of Viscous Incompressible Flow.} 2nd english
edition, Gordon and Breach, New York, 1969.

\bibitem{lemarie02} P. G. Lemari\'{e}-Rieusset. \textit{Recent developments in the Navier-Stokes problem.} 
Chapman \& Hall, Boca Raton, 2002.

\bibitem{MS}  P. Maremonti; V. A. Solonnikov. {\it On nonstationary Stokes problem in exterior domains.} Ann. Scuola Norm. Sup. Pisa Cl. Sci. (4) 24 (1997), no. 3, 395--449. 

\bibitem{majda93} A. Majda. {\it Remarks on weak solutions for vortex sheets with a distinguished sign.} Indiana Univ. Math. J. {\bf 42} (1993) 921--939.

\bibitem{MB02} A. Majda; A. Bertozzi. \textit{Vorticity and incompressible flow.} Cambridge Texts in Applied Mathematics 27, Cambridge University Press, Cambridge, 2002.

\bibitem{MNP00} V. Mazy'a, S. Nazarov; B. Plamenevskij. \textit{Asymptotic theory of elliptic boundary 
value problems in singularly perturbed domains, Vols 1 and 2.} Birkhauser Verlag, Basel, 2000.

\bibitem{saffman} P. Saffman. \textit{Vortex dynamics.} Cambridge Monographs on Mechanics and Applied Mathematics, 
Cambridge University Press, New York, 1992. 

\bibitem{serrin62} J. Serrin. {\it On the interior regularity of weak solutions of the Navier-Stokes equations.}
Arch. Rat. Mech. Anal. {\bf 9} (1962) 187--195.

\bibitem{temam84} R. Temam. {\it Navier-Stokes equations. Theory and numerical analysis. With an appendix by F. Thomasset.} Third edition, Studied in Mathematics and its Applications, 2, North-Holland Publishing Co., Amsterdam, 1984.   

\end{thebibliography}
\end{document}